\documentclass[11pt,reqno]{amsart}
\usepackage[french,english]{babel}
\usepackage[utf8]{inputenc}
\usepackage[T1]{fontenc}            
\usepackage{amsmath}
\usepackage{graphicx}
\usepackage{amssymb}
\usepackage{amsthm}
\usepackage{comment}
\usepackage{pgf, tikz}
\usepackage{amsfonts,mathrsfs}
\usepackage{thmtools}
\usepackage{geometry}
\usepackage{stmaryrd}
\usepackage{xcolor}
\usepackage{bm}
\usepackage{hyperref}
\hypersetup{colorlinks = true , linkcolor=blue, citecolor=blue}
\numberwithin{equation}{section}

\theoremstyle{definition}

\declaretheorem[numberwithin=section,title=Definition]{Def}
\declaretheorem[sibling=Def,title=Proposition]{Prop}
\declaretheorem[sibling=Def,title=Corollary]{Cor}
\declaretheorem[sibling=Def,title=Theorem]{Thm}
\declaretheorem[sibling=Def,title=Lemma]{Lemme}

\declaretheorem[sibling=Def,title=Remark]{Rq}
\declaretheorem[sibling=Def,title=Example]{Ex}

\newenvironment{dem}{\vspace{6pt} \paragraph{\textbf{Proof.} \hspace{3pt}}}{\hfill$\square$\newline }

\renewcommand{\P}{\mathbf{P}}
\newcommand{\PP}{\mathscr{P}}
\newcommand{\PPP}{\mathcal{P}}
\newcommand{\FF}{\mathcal{F}}
\newcommand{\E}{\mathbf{E}}

\newcommand{\LL}{\mathcal{L}}
\newcommand{\KK}{\mathcal{K}}
\newcommand{\R}{\ensuremath{\mathbf{R}}\xspace}

\newcommand{\N}{\ensuremath{\mathbf{N}}\xspace}

\newcommand{\CC}{\mathcal{C}}
\newcommand{\BB}{\mathcal{B}}

\newcommand{\1}{\mathbf{1}}

\newcommand{\del}{\frac{\delta u}{\delta  m}}

\begin{document}

	\title{Itô-Krylov's formula for a flow of measures}

\author{Thomas Cavallazzi}
	\address{Université de Rennes 1, CNRS, IRMAR - UMR 6625, F-35000 Rennes, France}
	\email{thomas.cavallazzi@univ-rennes1.fr}
	\date{November 07, 2022}

	\begin{abstract}
		We prove Itô's formula for the flow of measures associated with an Itô process having a bounded drift and a uniformly elliptic and bounded diffusion matrix, and for functions in an appropriate Sobolev-type space. This formula is the almost analogue, in the measure-dependent case, of the Itô-Krylov formula for functions in a Sobolev space on $\R^+ \times \R^d$. 
	\end{abstract}
 \maketitle

	\section{Introduction}\
	
	We fix $(\Omega,\FF,(\FF_t)_{t \geq 0},\P)$ a filtered probability space satisfying the usual conditions. Let $T>0$ be a finite horizon of time, $d,d_1 \in \N^*$ with $d_1 \geq d,$ and $(B_t)_{t \geq 0} $ a $(\FF_t)_{t \geq 0}$-Brownian motion of dimension $d_1$. We consider the Itô process on $\R^d$ defined, for $t \in [0,T],$ by  \begin{equation}\label{itoprocess}
	X_t := X_0 + \int_0^t b_s \, ds + \int_0^t \sigma_s \, dB_s,\end{equation} where $X_0 \in L^2(\Omega,\FF_0;\R^d)$, $b: [0,T]\times \Omega\rightarrow \R^d$ and $\sigma : [0,T]\times \Omega \rightarrow \R^{d \times d_1}$ are progressively measurable processes. In the following, we will denote by $\mu_t$ the law of $X_t$ and by $a$ the matrix $\sigma \sigma^*.$\\
	
	Let us fix a real-valued function $u$ defined on the 2-Wasserstein space $\PPP_2(\R^d),$ i.e. the space of probability measures on $\R^d$ having a finite moment of order $2.$ In this paper, we are interested in Itô's formula for $u$ and the flow of probability measures $(\mu_t)_{t \in [0,T]}.$ This formula describes the dynamics of $t \mapsto u(\mu_t),$ essentially by computing its derivative (see \eqref{formulaito11} below). It has a wide range of applications for example in Mean-Field Games, McKean-Vlasov's control problems, McKean-Vlasov Stochastic Differential Equations (SDEs) but also in the study of interacting particle systems and the propagation of chaos. These applications will be detailed below.\\
	
	 Itô's formula for a flow of measures naturally requires differential calculus on the space of measures $\PPP_2(\R^d).$ We will use the linear (functional) derivative, which is a standard notion of differentiability for functions of measures relying on the convexity of $\PPP_2(\R^d).$ The function $u$ admits a linear derivative if there exists a real-valued and continuous function $\del$ defined on $\PPP_2(\R^d)\times \R^d,$ at most of quadratic growth with respect to the space variable uniformly on each compact set of $\PPP_2(\R^d),$ and such that for all $\mu,\nu \in \PPP_2(\R^d)$ $$ u(\mu) - u(\nu) = \int_0^1 \int_{\R^d} \del(t\mu + ( 1-t) \nu)(v) \, d (\mu - \nu)(v) \, dt.$$

	The standard Itô formula for a flow of measures can be found in \cite{BuckdahnMeanfieldstochasticdifferential2014} (see Theorem 6.1) or in Section $3$ of \cite{ChassagneuxProbabilisticapproachclassical2015} and Chapter 5 of \cite{CarmonaProbabilisticTheoryMean2018} (see Theorem 5.99) under less restrictive assumptions. It states that for all $t \in [0,T]$  \begin{equation} 
	\label{formulaito11} u(\mu_t) = u(\mu_0) + \int_0^t \E \left( \partial_v \del (\mu_s)(X_s)\cdot b_s\right) \,ds + \frac{1}{2}  \int_0^t \E \left( \partial^2_v \del (\mu_s)(X_s)\cdot a_s\right) \,ds,  \end{equation}
	where $x\cdot y$ denotes the usual scalar product of two vectors $x,y \in \R^d$ and $A\cdot B := \text{Tr}(A^*B)$ the usual scalar product of two matrices $A,B \in \R^{d\times d}.$ The common point between these results is that the function $u$ has to be $\CC^2$ in some sense. More precisely, it is always assumed that for all $\mu \in \PPP_2(\R^d),$ the linear derivative $\del (\mu)(\cdot)$ belongs to $\CC^2(\R^d)$  or equivalently that the L-derivative $\partial_{\mu}u(\mu)(\cdot)$ belongs to $\CC^1(\R^d)$ (see below for the definition of the L-derivative and its link with the linear derivative). This paper aims at proving Itô's formula \eqref{formulaito11} for functions $u$ having a linear derivative $\del$ that is not $\CC^2$ with respect to the space variable.\\
	
We now fix the assumptions on the Itô process $(X_t)_{t \in [0,T]}.$ In this paper, we always assume that the drift $b$ and the diffusion matrix $\sigma$ in \eqref{itoprocess} satisfy the following properties. \begin{enumerate}
		\item[] \textbf{(A)}  There exists $K>0$ such that almost surely  $$  \forall t \in [0,T], \, |b_t| + | \sigma_t| \leq K.$$
		\item[] \textbf{(B)} There exists $\delta >0$ such that almost surely  $$ \forall t \in [0,T], \, \forall \lambda \in \R^d, \,  a_t\lambda\cdot\lambda  \geq \delta |\lambda|^2.$$ 
	\end{enumerate}
\vspace{4pt}

	  Assumptions \textbf{(A)} and \textbf{(B)} stem from Section $2.10$ of \cite{KrylovControlledDiffusionProcesses2009}. Therein, Krylov deals with controlled diffusion processes and needs to apply the standard Itô formula for the so-called pay-off function which is not $\CC^2.$ That is why he proves an extension of the classical Itô formula for the Itô process $(X_t)_{t \in [0,T]}$ satisfying Assumptions \textbf{(A)} and \textbf{(B)}, and for a function $g: \R^d \rightarrow \R$ belonging to an appropriate Sobolev space. The crucial point is that $(X_t)_t$ satisfies the non-degeneracy Assumption \textbf{(B)}. It ensures that the noise does not degenerate and allows to produce a regularizing effect. Let us explain how. The non-degeneracy assumption leads to Krylov's inequality (see Theorem \ref{inegalitekrylov} taken from Section $2.3$ of \cite{KrylovControlledDiffusionProcesses2009}). This inequality, in turn, implies that for almost all $t \in [0,T],$ $\mu_t,$ the law of $X_t,$ has a density $p(t,\cdot)$ with respect to the Lebesgue measure (see Proposition \ref{densiteexistence}). Moreover, this density belongs to $L^{(d+1)'}([0,T]\times \R^d),$ where $(d+1)'$ denotes the conjugate exponent of $d+1$ defined in Section \ref{notations}. The existence of densities together with the integrability property permit to assume Sobolev regularity for the function $g.$ More precisely, Itô-Krylov's formula is established under the assumption that $g$ is continuous on $\R^d$ and that $\nabla g$ belongs to the Sobolev space $W_{\text{loc}}^{1,k}(\R^d),$ for $ k \geq d+1,$ i.e. that $\nabla g$ and $\nabla^2 g$ are in $L_{\text{loc}}^{k}(\R^d)$ (see Section 2.10 of \cite{KrylovControlledDiffusionProcesses2009}).   \\

	  Our goal here is to take advantage of the regularizing effect of the noise, stemming from the existence of the densities $p(t,\cdot)$ and their integrability property, to establish an analogue of Itô-Krylov's formula in the measure-dependent case. Looking at Itô's formula for a flow of measures \eqref{formulaito11}, the regularizing effect comes from the presence of expectations which average, with respect to the space variable, the derivatives of $\del$ on all the trajectories of $(X_t)_t.$ Indeed, the regularization by noise will only appear through the space variable of the linear derivative but not through its measure variable. This is not surprising since the space of measures $\PPP_2(\R^d)$ is somehow infinite dimensional while the noise is of finite dimension. Thus, we cannot expect a true regularization in the measure variable of $\del.$ The fact that a finite dimensional noise cannot have a complete regularizing effect in the space $\PPP_2(\R^d)$ is explained in \cite{marx2020infinitedimensional} in the context of McKean-Vlasov SDEs.  \\
	  
	   In order to prove Itô's formula \eqref{formulaito11} for $u$, it is clear that $u$ needs to admit a linear derivative with at least distributional derivatives of order $1$ and $2$ with respect to the space variable in $L^k(\R^d)$ for some $k,$ as for the standard Itô-Krylov formula. Let us describe more precisely our assumptions on $u.$ As said before, for almost all $t \in [0,T],$ the law $\mu_t$ has a density $p(t,\cdot)$ such that $p$ belongs to $L^{(d+1)'}([0,T]\times \R^d).$ Denoting by $\PP(\R^d)$ the space of measures $\mu \in \PPP_2(\R^d)$ having a density with respect to the Lebesgue measure in $ L^{(d+1)'}(\R^d),$ our assumptions on the derivatives of $\del(\mu)(\cdot)$ are only made for measures $\mu$ belonging to $\PP(\R^d).$ This is natural since for almost all $t \in [0,T],$ $\mu_t$ belongs to $\PP(\R^d),$ and the derivatives of $\del$ are evaluated along the flow $(\mu_t)_{t \in [0,T]}$ and integrated in time. Moreover, because of the integrability property of the densities $p(t,\cdot),$ the derivatives of $\del (\mu)(\cdot)$ do not need to be defined and continuous on the whole space $\R^d$ because they are somehow integrated against the densities $p(t,\cdot)$ (see \eqref{formulaito11}). We say "somehow" because it is not completely the case since $b$ and $a$ are random. But as they are bounded, we can omit them in some sense. More precisely, the integrability property of the densities leads us to assume that $u$ admits a linear derivative such that for all $\mu \in \PP(\R^d),$ $\partial_v \del(\mu)(\cdot)$ belongs to the Sobolev space $W^{1,k}(\R^d)$ defined in Section \ref{notations}, with $k \geq d+1.$ This is exactly the same condition as in the standard Itô-Krylov formula, except that we replace $W_{\text{loc}}^{1,k}(\R^d)$ by $W^{1,k}(\R^d).$ This is essentially explained by the expectations in Itô's formula \eqref{formulaito11}. Indeed, the process $(X_t)_t$ cannot be localized by stopping times. Moreover, we assume that the map $\mu \in \PP(\R^d) \mapsto \partial_v \del(\mu)(\cdot) \in W^{1,k}(\R^d)$ is continuous for a distance on $\PP(\R^d)$ satisfying the assumptions of Definition \ref{defespaceP}. This continuity assumption can be interpreted as the fact that the noise has no regularizing effect in the measure variable of the linear derivative, as explained above. The precise assumptions of our Itô-Krylov's formula are given in Definition \ref{sobolevW1} and Theorem \ref{ito_krylov}. Eventually, we extend  in Theorem \ref{extensionito} our formula to functions depending also on the time and space variables satisfying the assumptions of Definition \ref{sobolevW}.\\

	 	We now focus on some applications of Itô's formula for a flow of measures. This one has been developed with the increasing interest for Mean-Field Games and McKean-Vlasov SDEs over the last decade. Mean-Field Games were initiated independently by Caines, Huang and Malhame in \cite{CainesLargepopulationstochastic2006} and  by Lasry and Lions in \cite{LasryMeanfieldgames2007}. The notion of Master equations has been introduced by Lions in his lectures at Collège de France \cite{LionscourscollegedeFrance} in order to describe Mean-Field Games. Master equations are Partial Differential Equations (PDEs) on the space of probability measures and can be derived with the help of Itô's formula. We refer to Lions' lectures \cite{LionscourscollegedeFrance}, the notes written by Cardialaguet \cite{CardaliaguetnotesMFG2013}, and the books of Carmona and Delarue \cite{CarmonaProbabilisticTheoryMean2018,CarmonaProbabilisticTheoryMean2018a} for more details on Mean-Field Games and Master equations. We also mention Bensoussan, Frehse and Yam \cite{BensoussanMasterEquationMean2014} and Carmona, Delarue \cite{CarmonaMasterEquationLarge2014} where Master equations are derived, with the help of Itô's formula in \cite{CarmonaMasterEquationLarge2014}. The question of existence and uniqueness of classical solutions to Master equations was addressed by Cardaliaguet, Delarue, Lasry and Lions in \cite{CardaliaguetMasterEquationConvergence2015} and by Chassagneux, Crisan and Delarue in \cite{ChassagneuxProbabilisticapproachclassical2015}. From a different point of view, Mou and Zhang deal with  the well-posedness of Master equations in some weaker senses in \cite{MouWellposednessSecondOrder2020}. \\

  Moreover, Itô's formula appears to be the natural way to connect a McKean-Vlasov SDE (more precisely the associated semigroup $(P_t)_t$ acting on the space of functions of measures) to a PDE on the space of probability measures (the Master equation) in the same manner as for classical SDEs. It turns out to be a crucial tool to study the stochastic flow generated by a McKean-Vlasov SDE, as explained in Chapter $5$ of \cite{CarmonaProbabilisticTheoryMean2018}. The link between McKean-Vlasov SDEs and PDEs on the space of measures is at the heart of the work of Buckdahn, Li, Peng and Rainer \cite{BuckdahnMeanfieldstochasticdifferential2014}  where the authors prove that the PDE admits a unique classical solution expressed with the flow of measures associated with the McKean-Vlasov SDE. Moreover, in the parallel work \cite{ChassagneuxProbabilisticapproachclassical2015}, Chassagneux, Crisan and Delarue adopt a similar approach and study the flow generated by a forward-backward stochastic system of McKean-Vlasov type under weaker assumptions on the coefficients of the equation. Both works are motivated by Mean-Field Games, and Itô's formula plays a key role. In \cite{CrisanSmoothingpropertiesMcKeanVlasov2017}, Crisan and McMurray prove that the Master equation admits a unique classical solution for some irregular terminal condition using Malliavin calculus. They point out a smoothing effect concerning the differentiability of the solution with respect to the measure even though there is no noise in the measure direction. Furthermore, the problem of propagation of chaos for the interacting particles system associated with the McKean-Vlasov SDE can also be addressed with the help of the associated PDE on the space of measures (see Chapter 5 of \cite{CarmonaProbabilisticTheoryMean2018}). It allows to obtain quantitative weak propagation of chaos estimates between the law of the solution to the McKean-Vlasov SDE and the empirical measure of the associated particle system. This approach was adopted for example by Chaudru de Raynal and Frikha in \cite{deraynal2021wellposedness,frikha2021backward}, by Delarue and Tse in \cite{delarue2021uniform} and by Chassagneux, Szpruch and Tse in \cite{chassagneux2019weak}. Let us also mention that the Master equation satisfied by the semigroup has been recently used by Jourdain and Tse in \cite{JourdainCentrallimittheorem2021} to study the mean-field fluctuation (CLT) of an interacting particle system. Finally, Itô's formula for a flow of measures is also important to deal with McKean-Vlasov control problems because it allows to derive a dynamic programming principle describing the value function of the problem as presented in Chapter $6$ of \cite{CarmonaProbabilisticTheoryMean2018}.\\

	 Recently, Itô's formula has been extended to flows of measures generated by càdlàg semi-martingales. It was achieved independently by Guo, Pham and Wei in \cite{GuoItoFormulaFlow2020}, who studied McKean-Vlasov control problems with jumps and by Talbi, Touzi and Zhang in \cite{talbi2021dynamic} who worked on mean-field optimal stopping problems. In both works, dynamic programming principles are established thanks to Itô's formula for a flow of measures. Finally, we also mention that several Itô-Wentzell-Lions formulae for functional random fields of Itô type depending on measure flows have been established by dos Reis and Platonov in \cite{dosReisItoWentzellLionsFormulaMeasure2022}. \\
	 
	 Let us explain our choice to work with the linear derivative. Indeed, the L-derivative, which was introduced by Lions in his lectures at Collège de France \cite{LionscourscollegedeFrance}, is also well-adapted to establish Itô's formula for a flow of measures. We say that $u$ is L-differentiable if its lifting defined by $$ \tilde{u} : X \in L^2(\Omega;\R^d) \mapsto u(\mathcal{L}(X)) \in \R,$$ where $ \mathcal{L}(X)$ denotes the law of $X,$ is Fréchet differentiable on $L^2(\Omega;\R^d).$ Moreover, there exists a $\R^d$-valued function $\partial_{\mu} u$ defined on $\PPP_2(\R^d) \times \R^d$ such that the gradient of $\tilde{u}$ at $X \in L^2(\Omega;\R^d)$ is given by the random variable $\partial_{\mu}u (\mathcal{L}(X))(X).$ The function $\partial_{\mu}u$ is called the L-derivative of $u.$ The advantage of the L-derivative is that it permits to use standard tools of differential calculus on Banach spaces. Of course, there is a link between the L-derivative and the linear derivative of $u.$ Indeed, in general, the L-derivative $\partial_{\mu} u(\mu)(\cdot)$ is equal to the gradient of the linear derivative  $\partial_v \del (\mu)(\cdot)$ (see Propositions 5.48 and 5.51 in \cite{CarmonaProbabilisticTheoryMean2018} for the precise assumptions). Under our assumptions presented above, Sobolev embedding theorem ensures that for all $\mu \in \PP(\R^d),$ $\del(\mu)(\cdot)$ belongs to $\CC^1(\R^d;\R),$ and that $\partial_v \del(\mu)(\cdot)$ is continuous and bounded on $\R^d.$ We would be tempted to deduce that $u$ admits a L-derivative given, as recalled above, by  $\partial_v \del(\mu)(\cdot).$ However, this term is assumed to exist only for measures $\mu \in \PP(\R^d)$ and not for $\mu \in \PPP_2(\R^d).$ This is the case in Example \ref{exquadratic}, where this term is not well-defined for any $\mu \in \PPP_2(\R^d)$ (see Remark \ref{rqchoicelinearderivative}). It seems therefore more restrictive to work with the L-derivative and thus justifies our choice to work with the linear derivative. \\

		   The paper is organized as follows. Section \ref{notations} gathers some notations and definitions used throughout the paper.  In Section \ref{sectionspacesandformula}, more precisely in Definitions \ref{sobolevW1} and \ref{sobolevW}, we define the spaces of functions for which we will establish Itô-Krylov's formula. These formulas are given in Theorem \ref{ito_krylov} for functions defined on $\PPP_2(\R^d)$ and in Theorem \ref{extensionito} for functions depending also on the time and space variables. Moreover, we give examples of functions for which our formulas hold and we discuss our assumptions through them. The proofs of these examples are postponed to Appendix \ref{sectionappendix} for ease of reading. In Section \ref{sectionpreliminaries}, we give some preliminary results. We start with Krylov's inequality and its consequences on the existence of densities for the flow of measures $(\mu_t)_{t\in[0,T]}$ in Proposition \ref{densiteexistence}. Then we recall some classical results on convolution and regularization. Finally, Sections \ref{sectionproof1} and \ref{sectionproof2} are respectively dedicated to the proofs Theorems \ref{ito_krylov} and \ref{extensionito}.\\

	\section{Notations and definitions}\label{notations}
	
	\subsection{General notations}
	
		\noindent Let us introduce some notations used several times in the article.
	\begin{enumerate}
		
		\item[-]$B_R$ is the open ball centered at $0$ and of radius $R$ in $\R^d$ for the euclidean norm. 
		
		\item[-]  $p'$ is the conjugate exponent of $p \in [1,+\infty]$, defined by $\frac1p + \frac{1}{p'} = 1.$
		
		\item[-] $L^p_{\text{loc}}(\R^d)$ is the space of functions $f$ such that for all $R>0$, $f\in L^p(B_R).$	\item[-] $W^{m,k}(\mathcal{O})$ is the Sobolev space of functions $u \in L^{k}(\mathcal{O})$ admitting distributional derivatives of order between $1$ and $m$ in $ L^{k}(\mathcal{O}),$ where $\mathcal{O}$ is open in $\R^d.$ It is equipped with the norm  $$ \|u\|_{W^{m,k}(\mathcal{O})} =  \sum_{\alpha \in \N^d, \, |\alpha| \leq m} \|\partial^{\alpha} u \|_{L^k(\mathcal{O})}.$$
		
		\item[-]$W^{m,k}_{\text{loc}}(\R^d)$ is the space of functions $u$ such that for all $R>0$, $u$ belongs to $W^{m,k}(B_R).$
		
		\item[-]  $(\rho_n)_n$ is a mollifying sequence on $\R^d$, that is a sequence of non-negative $\CC^{\infty}$ functions, such that for all $n$, $\int_{\R^d} \rho_n(x) \, dx =1$ and $\rho_n$ is equal to $0$ outside $B_{1/n}.$ We assume that $\rho_n(x) = \rho_n(-x)$ for all $x.$
		
		\item[-] $*$ denotes the convolution of two functions, when it is well-defined, or two probability measures.
		
		\item[-] $\BB(E)$ is the Borel $\sigma$-algebra where $E$ is a metric space. 
		\item[-] $A^*$ denotes the transpose of the matrix $A \in \R^{d\times d}.$
		\item[-] $A\cdot B$ denotes the usual scalar product of two matrices $A,B \in \R^{d \times d}$ given by $A\cdot B := \text{Tr}(A^*B).$
		\item[-] $\PP(\R^d)$ is defined in Definition \ref{defespaceP}.
		\item[-] $\mathcal{W}_1(\R^d)$ is defined in Definition \ref{sobolevW1}.
		\item[-] $\mathcal{W}_2(\R^d)$ is defined in Definition \ref{sobolevW}.
		
	\end{enumerate}

\subsection{Spaces of measures and linear derivative}
	
	The set $\PPP(\R^d)$ is the space of probability measures on $\R^d$ equipped with the topology of weak convergence. The Wasserstein space $\PPP_2(\R^d)$ denotes the set of measures $\mu \in \PPP(\R^d)$ such that $\int_{\R^d} |x|^2 \,d\mu(x) < + \infty,$ equipped with the $2$-Wasserstein distance $W_2$ defined for $\mu, \nu \in \PPP_2(\R^d)$ by $$ W_2(\mu,\nu) = \inf_{\pi \in \Pi(\mu,\nu)} \left(\int_{\R^d\times \R^d} |x-y|^2 \, d\pi(x,y)\right)^{1/2},$$ where $\Pi(\mu,\nu)$ is the subset of $\PPP_2(\R^d \times \R^d)$ with marginal distributions $\mu$ and $\nu.$ We will work with the standard notion of linear derivative for functions of measures.
	
	\begin{Def}[Linear derivative] A function $ u: \PPP_2(\R^d) \rightarrow \R$ is said to have a linear derivative if there exists a continuous function $(\mu,v) \in \PPP_2(\R^d) \times \R^d \mapsto \del(\mu)(v) \in \R,$ satisfying the following properties. \begin{enumerate}
			\item For all compact $ \KK \subset \PPP_2(\R^d)$  $  \displaystyle\sup_{v \in \R^d}\displaystyle\sup_{\mu \in \KK} \left\{ (1 + |v|^2)^{-1}\left| \del (\mu)(v)\right|\right\} < + \infty.$\item For all $\mu,\nu \in \PPP_2(\R^d),$ $u(\mu) - u(\nu) = \displaystyle\int_0^1 \int_{\R^d} \del(t\mu + (1-t) \nu)(v) \, d(\mu - \nu)(v) \, dt.$ 
		\end{enumerate}
	\end{Def}
	
	\begin{Rq}\label{Rqlinearderivative}Instead of the second point of the previous definition, it is equivalent to assume that for all $\mu,\nu \in \PPP_2(\R^d),$ $t \in [0,1] \mapsto u(t\mu + (1-t)\nu)$ is of class $\CC^1$ with $$\forall t \in [0,1], \, \frac{d}{dt} u(t\mu + (1-t)\nu) = \int_{\R^d} \del(t\mu + (1-t)\nu)(v) \, d(\mu-\nu)(v).$$ 
	\end{Rq}
	One can find more details in Chapter $5$ of \cite{CarmonaProbabilisticTheoryMean2018}, in particular the connection with the $L$-derivative.\\
	
	 Let us fix $(\rho_n)_n$ a mollifying sequence on $\R^d$, that is a sequence of non-negative $\CC^{\infty}$ functions, such that for all $n$, $\int_{\R^d} \rho_n(x) \, dx =1$ and $\rho_n$ is equal to $0$ outside $B_{1/n}.$ We assume that $\rho_n(x) = \rho_n(-x)$ for all $x.$
	
	\begin{Def}\label{defespaceP}
		Let us define $\PP(\R^d)$ as the space of measures $\mu \in \PPP_2(\R^d)$ which admit a density $\frac{d\mu}{dx}$ with respect to the Lebesgue measure belonging to $L^{(d+1)'}(\R^d).$ We endow $\PP(\R^d)$ with a general distance $d_{\PP}$ satisfying the following properties. \begin{enumerate}
			\item[] \textbf{(H1)} For any $n \geq 1$, $\mu \in (\PPP_2(\R^d),W_2) \mapsto \mu * \rho_n \in (\PP(\R^d),d_{\PP})$ is continuous.
			
			\item[] \textbf{(H2)} For any $\mu \in \PP(\R^d),$ $\mu * \rho_n \underset{n \rightarrow + \infty}{\longrightarrow} \mu$ for $d_{\PP}.$ 
		\end{enumerate}
	\end{Def}
	
	Note that for all $n \geq 1$ and for all $\mu \in \PPP_2(\R^d)$, $\mu * \rho_n \in \PP(\R^d).$ Indeed, its density is given by $x \mapsto \rho_n * \mu(x) = \int_{\R^d} \rho_n(x-y) \, d\mu(y).$ Jensen's inequality ensures that it belongs to $L^{(d+1)'}(\R^d).$ Considering the space $(\PP(\R^d),d_{\PP})$ comes in a natural way with Assumptions \textbf{(A)} and \textbf{(B)} on the Itô process $X.$ As explained in the introduction, it implies the existence of a density $p \in L^{1}([0,T] \times \R^d;\R^+) \cap L^{(d+1)'}([0,T] \times \R^d;\R^+)$ such that for almost all $t \in [0,T],$ the law of $X_t$ is equal to $ p(t,\cdot)\,dx$ and belongs to $\PP(\R^d)$ (see Proposition \ref{densiteexistence}). Let us give two examples for the distance $d_{\PP}.$ 
	
	\begin{Ex}\label{choicedistance} The Wasserstein distance $W_2$ clearly satisfies Assumptions \textbf{(H1)} and \textbf{(H2)} in Definition \ref{defespaceP}. Another family of examples is given by the distance $d_k$ defined, for $ k \in [d+1,+\infty[, \mu,\nu \in \PP(\R^d),$ by
		$$ d_k(\mu,\nu) = \left\Vert \frac{d\mu}{dx} - \frac{d\nu}{dx} \right\Vert_{L^{k'}(\R^d)}.$$ 
		
	\end{Ex}
	Note that $d_k$ is well-defined since for any $\mu \in \PP(\R^d)$, $\frac{d\mu}{dx} \in L^1(\R^d) \cap L^{(d+1)'}(\R^d)$ which is included in $L^{k'}(\R^d)$ by interpolation. The proof is postponed to the Appendix (Section \ref{proofchoicedistance}).\\
	
	\section{Itô-Krylov's formula, ah-hoc spaces of functions and examples}\label{sectionspacesandformula}

 Let us introduce now the Sobolev-type space of functions on $\PPP_2(\R^d)$ for which we will prove Itô's formula for a flow of measures.
	
	\begin{Def}\label{sobolevW1}
			
		Let $\mathcal{W}_1(\R^d)$ be the space of continuous functions $u: \PPP_2(\R^d) \rightarrow \R$ having a linear derivative $\del$ such that for all $\mu \in \PP(\R^d)$, the function $ \del (\mu)(\cdot)$ admits distributional derivatives of order $1$ and $2$ in $L^{k}(\R^d)$, for a certain $k \geq d+1,$ and satisfies the following properties. \begin{enumerate}
		
			\item  The map $  \mu \in (\PP(\R^d),d_{\PP}) \mapsto \partial_v \del (\mu)(\cdot) \in \left(W^{1,k}(\R^d)\right)^d$  is continuous  for a certain distance $d_{\PP}$ satisfying \textbf{(H1)} and \textbf{(H2)}.
			\item There exists $\alpha \in \N$ such that $ k\geq (1+\alpha)d$ and for all compact $\KK \subset \PPP_2(\R^d)$ and for any $ \mu \in \KK \cap \PP(\R^d)$ $$  \left\Vert \partial_v \del (\mu)(\cdot) \right\Vert_{L^k(\R^d)} + \left\Vert \partial_v^2 \del (\mu)(\cdot) \right\Vert_{L^k(\R^d)} \leq C_{\KK} \left( 1 + \left\Vert \frac{d\mu}{dx} \right\Vert_{L^{k'}(\R^d)}^{\alpha}\right).$$
			\end{enumerate}
	\end{Def}

	\begin{Rq}\label{rqW1}

-The space $\mathcal{W}_1(\R^d)$ contains the functions which satisfy Assumption $(1)$ in Definition \ref{sobolevW1} with   $(\PPP_2(\R^d),W_2)$ instead of $(\PP(\R^d),d_{\PP}).$ Indeed, the second point is clearly satisfied with $\alpha =0$ since $\KK$ is compact. \\

\noindent-Assumption $(2)$ in Definition \ref{sobolevW1} allows to control the growth of $\left\Vert \partial_v \del(\mu)(\cdot)\right\Vert_{W^{1,k}(\R^d)}$ with respect to the measure $\mu.$ It allows us to take advantage of the continuity of the flow in $\PPP_2(\R^d)$ (because the control is assumed on compact subsets of $\PPP_2(\R^d)$), but also of its integrability properties proved in Lemmas \ref{integrability1} and \ref{integrability2}. The form of the inequality suggests the integration of functions in $L^k(\R^d)$ with respect to $\mu,$ at least when the function $u$ is linear in $\mu.$\\

\noindent-Sobolev embedding theorem (see Corollary 9.14 in \cite{BrezisFunctionalAnalysisSobolev2010}) ensures that for all $\mu \in \PP(\R^d)$, $\del (\mu)(\cdot)$ belongs to $\CC^1(\R^d;\R)$ and that $\partial_v \del(\mu)(\cdot)$ is bounded and $\gamma$-Hölder, where $\gamma:= 1 - \frac{d}{k}$. Note that we do not need that $\del(\mu)(\cdot) \in W^{2,k}(\R^d)$ since there is no integrability assumption made on the linear derivative.

\end{Rq} 

Having this definition at hand, we can now state Itô-Krylov's formula for functions in $\mathcal{W}_1(\R^d).$ 

\begin{Thm}[Itô-Krylov's formula]\label{ito_krylov}\
	Let $ u $ be a function in $\mathcal{W}_1(\R^d),$  which was defined in Definition \ref{sobolevW1}. We have for all $t \in [0,T]$  \begin{align}\label{formulaito1} u(\mu_t) &= u(\mu_0) + \int_0^t \E \left( \partial_v \del (\mu_s)(X_s)\cdot b_s\right) \,ds + \frac{1}{2}  \int_0^t \E \left( \partial^2_v \del (\mu_s)(X_s)\cdot a_s\right) \,ds,
	\end{align}
	where  $\partial^2_v \del (\mu_s)(X_s)\cdot a_s := \text{Tr}\Big(\partial^2_v \del (\mu_s)(X_s)a_s\Big)$ is the usual scalar product on $\R^{d\times d}.$ \\
	
\end{Thm}

\begin{Rq}
	Notice that a function  $u \in\mathcal{W}_1(\R^d)$ is assumed to have a linear derivative on the whole space $\PPP_2(\R^d)$. This seems a bit strong at first sight in comparison with the assumptions on its spatial derivatives that are only made for measures $\mu \in \PP(\R^d).$ Indeed, we could consider working with a linear derivative defined only on the space of densities, as done for example in \cite{BensoussanMasterEquationMean2014}. However, in order to establish Itô-Krylov's formula by regularization, the function $u$ needs to be continuous on the whole space $\PPP_2(\R^d)$ and not only on $\PP(\R^d)$. Indeed, the flow $s \in [0,T] \mapsto \mu_s \in \PPP_2(\R^d)$ is continuous but $\mu_t$ does not necessarily belong to $\PP(\R^d)$ for all $t \in [0,T]$. This is proved only for almost all $t$. Thus, as the function $u$ has to be continuous on $\PPP_2(\R^d)$, we have chosen to assume the existence of a linear derivative on $\PPP_2(\R^d)$ even though we could have only required it on the space of densities.
\end{Rq}

Now, we focus on examples of functions belonging to $\mathcal{W}_1(\R^d).$ Let us start with the linear case. 

\begin{Ex}[Linear functional]\label{exlinear}
	Fix $g \in \CC^0(\R^d;\R)$ admitting a distributional derivative such that $\nabla g \in (W^{1,k}(\R^d))^d$ for some $k \geq d+1.$ Then, the function $$u : \left\{ \begin{array}{rll}
	\PPP_2(\R^d) &\rightarrow \R \\ \mu&\mapsto \displaystyle\int_{\R^d } g(x) \, d\mu(x),  \end{array} \right. $$ belongs to the space $\mathcal{W}_1(\R^d).$
\end{Ex}

Indeed, Sobolev embedding theorem (see Corollary 9.14 in \cite{BrezisFunctionalAnalysisSobolev2010}) implies that $\nabla g \in L^{\infty}(\R^d)$ since $ k \geq d+1.$ Thus $g$ is at most of linear growth so that for all $\mu \in \PPP_2(\R^d),$ $\del(\mu) = g,$ which clearly satisfies Assumptions $(1)$ and $(2)$ (with $\alpha = 0$) in Definition \ref{sobolevW1}.\\

	Let us now focus on the multi-linear case.
		
		\begin{Ex}[Polynomials on the Wasserstein space]\label{exquadratic}
			Fix $N\geq 2$ and $g \in \CC^0((\R^d)^N ;\R)$ such that \begin{enumerate}
				\item[-] there exists $C>0$ such that for all $\bm{x}=(x_1,\dots,x_N) \in (\R^d)^N, \,|g(\bm{x})| \leq C(1+|x_1|^2+ \dots +|x_N|^2),$
				\item[-] the distributional derivative $\nabla g$ belongs to $(W^{1,k}((\R^{d})^N))^{Nd}$ for a certain $k \in [Nd, + \infty[.$
				
			\end{enumerate}
			Then, the function $$u : \left\{ \begin{array}{rll}
			\PPP_2(\R^d) &\rightarrow \R \\ \mu&\mapsto \displaystyle\int_{(\R^d)^N} g(x_1,\dots,x_N) \, d\mu(x_1) \dots \, d\mu(x_N),  \end{array} \right. $$ belongs to the space $\mathcal{W}_1(\R^d)$ for $d_{\PP}=d_k.$
		\end{Ex}
	
	The proof is postponed to the Appendix  (Section \ref{proofexquadratic}).
		
		\begin{Rq}\label{rqchoicelinearderivative}
		
		 	- In Definition \ref{sobolevW1}, the distributional derivatives of the linear derivative $\del(\mu)$ are not necessarily integrable  functions for all $\mu \in \PPP_2(\R^d).$ Of course, in Example \ref{exlinear}, it is the case for all $\mu \in \PPP_2(\R^d)$ as the linear derivative does not depend on the measure $\mu.$ However, in Example \ref{exquadratic} for $N=2$, the linear derivative is given by \begin{equation}\label{linderquadra}
					\del (\mu)(v) = \int_{\R^d} g(v,y) \, d\mu(y) + \int_{\R^d} g(y,v) \, d\mu(y).\end{equation}   Formally, the derivative with respect to $v$ of the first integral in \eqref{linderquadra} is $$ \int_{\R^d} \partial_v g(v,y) \, d\mu(y).$$ This term is not well-defined for general measures $\mu \in \PPP_2(\R^d)$ because we have only assumed that $\nabla g\in (W^{1,k}(\R^{2d}))^{2d}$ with $k \geq 2d.$ Indeed, for $k=2d,$ we just know by Sobolev embedding theorem that $\nabla g$ belongs to $(L^r(\R^{2d})^{2d}$ with $r \in [2d,+\infty[$ (see Corollary 9.11 in \cite{BrezisFunctionalAnalysisSobolev2010}). As we will see in the proof (Section \ref{proofexquadratic} of the Appendix), it is well-defined as an integrable function of $v$ if we restrict to measures $\mu \in \PP(\R^d).$ This also justifies why we have chosen to work with the linear derivative instead of the L-derivative. Indeed, the L-derivative of $u$ would be equal to the gradient of the linear derivative $\partial_v \del(\mu)(\cdot),$ which is not well-defined for all $\mu \in \PPP_2(\R^d).$ Thus, the function $u$ does not need to be L-differentiable in the usual sense in our setting.\\
		
		\noindent	- Our assumptions on the derivatives of $\del$ in Definition \ref{sobolevW1} deal with $\PP(\R^d)$ instead of the whole space $\PPP_2(\R^d)$ essentially because in Itô's formula \eqref{formulaito1}, these derivatives only appear under integrals along the flow $(\mu_s)_{s \in [0,T]},$ which belongs to $\PP(\R^d)$ for almost all $s\in [0,T].$ However, we assume that $u$ is continuous on $\PPP_2(\R^d)$ since the flow $s \in [0,T] \mapsto \mu_s \in \PPP_2(\R^d)$ is continuous but $\mu_t$ does not necessarily belong to $\PP(\R^d)$ for all $t \in [0,T]$ .

		\end{Rq}
	
	The next example focuses on the particular case of convolution which has to be treated differently than in Example \ref{exquadratic} with $N=2$ because of the structure of the convolution which mixes the two variables. 
		
		\begin{Ex}\label{exconvol}
			Let $f \in \CC^0(\R^d;\R)$ be a function such that the distributional derivative $\nabla f$ belongs to  $(W^{1,k+1}(\R^d))^d,$ for a certain $ k\geq d.$ Then, the function  $$u: \left\{ \begin{array}{rll}
			\PPP_2(\R^d) &\rightarrow \R \\ \mu&\mapsto \displaystyle\int_{\R^d} f * \mu \, d\mu,  \end{array} \right.$$ belongs to $\mathcal{W}_1(\R^d)$ for $d_{\PP}=W_2.$
		\end{Ex}

	Here, the particular structure of convolution enables us to  work on the whole space $\PPP_2(\R^d)$ instead of $\PP(\R^d),$ as explained in the first point of Remark \ref{rqW1}. The proof is postponed to the Appendix (Section \ref{proofexconvol}). \\

		Finally, we give a non-linear example of functions belonging to $\mathcal{W}_1(\R^d)$.
		
		\begin{Ex}\label{ex-non-lin}
			Let $F \in \CC^1(\R ; \R)$ and $g \in \CC^0(\R^d;\R)$ be such that the distributional derivative $\nabla g$ belongs to $(W^{1,k}(\R^d))^d$ for some $k\geq d+1.$ Then $$u: \left\{ \begin{array}{cll}
			 \PPP_2(\R^d) &\rightarrow \R \\\mu&\mapsto F\left(\int_{\R^d} g \, d\mu \right)  \end{array} \right.$$ belongs to $\mathcal{W}_1(\R^d)$ for $d_{\PP}=W_2.$  
		\end{Ex}
	
	The proof is again postponed in the Appendix (Section \ref{proofex-non-lin}).\\

	We now deal with the extension of Itô's formula for functions depending also on the time and space variables. First, we define the space of functions generalizing the space $\mathcal{W}_1(\R^d).$ 
		
	\begin{Def}\label{sobolevW}
		Let $\mathcal{W}_2(\R^d)$ be the set of continuous functions $u: [0,T] \times \R^d \times \PPP_2(\R^d)\rightarrow \R$ satisfying the following properties for a certain distance $d_{\PP}$ satisfying \textbf{(H1)} and \textbf{(H2)}.
		
		\begin{enumerate}
			\item For all $(x,\mu) \in \R^d \times\PPP_2(\R^d)$, $u(\cdot,x,\mu)\in \CC^{1}$ and $\partial_t u$ is continuous on $[0,T] \times \R^d \times \PPP_2(\R^d)$. 
			\item There exists $k_1 \geq d+1$ such that for all $(t,\mu) \in [0,T] \times \PP(\R^d)$, $u(t,\cdot,\mu) \in W^{2,k_1}_{\text{loc}}(\R^d)$ and for all $ t \in [0,T]$ and $R>0$ 
			$$ \mu \in (\PP(\R^d), d_{\PP})\mapsto \partial_x u(t,\cdot,\mu) \in \left(W^{1,k_1}(B_R)\right)^d,$$ is continuous and  $\partial_x u$ and $\partial_x^2 u$ are measurable with respect to $(t,x,\mu) \in [0,T]\times \R^d\times \PP(\R^d) .$
			\item For all $(t,x) \in [0,T] \times \R^d$, $u(t,x,\cdot)$ admits a linear derivative $\del (t,x,\cdot)(\cdot)$ which is continuous on $[0,T]\times \R^d \times \PPP_2(\R^d)\times \R^d,$ and such that for all $\KK \subset  \R^d \times\PPP_2(\R^d)$ compact and $t\in [0,T]$, there exists $C>0$ such that for all $v \in \R^d$  $$   \sup_{(x,\mu) \in \KK }\left|\del (t,x,\mu)(v) \right| \, dx \leq C (1+|v|^2).$$\item There exists $k_2 \geq 2d$ such that for all $(t,\mu) \in [0,T] \times \PP(\R^d),$ $\del (t, \cdot,\mu)(\cdot)$ admits distributional derivatives with respect to $v$ of order $1$ and $2$ such that for all $t$ and $R>0 $ $$\mu\in  (\PP(\R^d),d_{\PP}) \mapsto \left(\partial_v \del (t,\cdot,\mu)(\cdot)\,,\, \partial^2_v \del (t,\cdot,\mu)(\cdot)\right) \in (L^{k_2}(B_R\times \R^d))^d \times (L^{k_2}(B_R\times \R^d))^{d \times d},$$ is continuous and measurable with respect to $(t,x,\mu,v) \in [0,T]\times \R^d \times \PP(\R^d) \times \R^d$.
			\item There exists $\alpha_1, \alpha_2 \in \N$ with $k_1\geq (2\alpha_1+1)d$, $ k_2 \geq (\alpha_2 +2)d$ such that for all $\KK \subset \PPP_2(\R^d)$ compact and $R>0$, there exists $C_{\KK,R}>0$ such that for all $\mu \in \KK \cap \PP(\R^d)$  $$\left\{ \begin{aligned}
				  &\sup_{t\leq T} \left\{ \left\Vert \partial_x u (t,\cdot,\mu) \right\Vert_{L^{k_1}(B_R)} +  \left\Vert \partial_x^2 u (t,\cdot,\mu) \right\Vert_{L^{k_1}(B_R)}\right\} \leq C_{\KK,R}  \left( 1 + \left\Vert \frac{d\mu}{dx} \right\Vert_{L^{k_1'}(\R^d)}^{\alpha_1}\right) \\&\sup_{t\leq T} \left\{ \left\Vert \partial_v \del (t,\cdot,\mu)(\cdot) \right\Vert_{L^{k_2}(B_R\times\R^d)} +  \left\Vert \partial_v^2 \del (t,\cdot,\mu)(\cdot) \right\Vert_{L^{k_2}(B_R\times\R^d)}\right\} \leq C_{\KK,R}  \left( 1 + \left\Vert \frac{d\mu}{dx} \right\Vert_{L^{k_2'}(\R^d)}^{\alpha_2}\right).  \end{aligned}   \right.$$

		\end{enumerate}
	\end{Def}

\begin{Rq}\label{CorW_2}

- The space $\mathcal{W}_2(\R^d)$ contains the functions satisfying the four first assumptions of Definition \ref{sobolevW} with $(\PP(\R^d),d_{\PP})$ replaced by $ (\PPP_2(\R^d),W_2)$ and also assuming that the functions in Assumptions $(2)$ and $(4)$ are continuous with respect to $(t,\mu)\in[0,T]\times \PPP_2(\R^d)$. Indeed, Assumption $(5)$ is automatically satisfied with $\alpha_1=\alpha_2 = 0$ because $\KK$ is compact. \\

\noindent- The bound in Assumption $(3)$ is quite natural. If the supremum in this bound was taken only over a compact set of $\PPP_2(\R^d)$, it would be the definition of the linear derivative. But we also need to control $\del$ locally uniformly in the space variable $x \in \R^d$ because of our regularization procedure through a convolution both in the space and measure variables. Assumptions $(2)$, $(4)$ and $(5)$ are generalizations of those in Definition \ref{sobolevW1} adapted to the presence of the space and time variables. In Assumption $(5)$, the condition on $k_2$ and  $\alpha_2$ changes a bit compared to the analogous assumption in Definition \ref{sobolevW1}, essentially because it deals with functions on $\R^{2d}$ instead of $\R^d.$ Let us mention that Assumption $(5)$ in Definition \ref{sobolevW} can be replaced by the integrability properties \eqref{bound1} established in Step $1$ of the proof of the next theorem (see Section \ref{sectionproof2}). 

\end{Rq} 

The next theorem is the natural extension of the formula for functions in $\mathcal{W}_2(\R^d).$  Let $(\eta_s)_{s \in [0,T]} $ and  $(\gamma_s)_{s \in [0,T]} $ be two progressively measurable processes, taking values respectively in $\R^d$ and $\R^{d \times d_1}$ and satisfying Assumptions \textbf{(A)} and \textbf{(B)}. We set, for all $t \leq T$  $$\xi_t = \xi_0 + \int_0^t \eta_s \, ds + \int_0^t \gamma_s \, dB_s,$$ where $\xi_0$ is a $\FF_0$-measurable random variable with values in $\R^d.$

\begin{Thm}[Extension of Itô-Krylov's formula]\label{extensionito}
	
	Let $u$ be a function in $\mathcal{W}_2(\R^d),$ which was defined in Definition \ref{sobolevW}. We have almost surely, for all $t \in [0,T]$ 
	\begin{align}\label{formulaito2}
	\notag u(t,\xi_t,\mu_t) &= u(0,\xi_0, \mu_0) + \int_0^t ( \partial_t u(s,\xi_s,\mu_s) + \partial_x u(s,\xi_s,\mu_s)\cdot\eta_s) \, ds + \frac12 \int_0^t \partial^2_x u(s,\xi_s,\mu_s)\cdot \gamma_s\gamma_s^* \, ds \\ &\quad+ \int_0^t \tilde{\E} \left(\partial_v \del (s,\xi_s,\mu_s)(\tilde{X_s})\cdot\tilde{b_s}\right) \, ds + \frac12 \int_0^t \tilde{\E} \left(\partial^2_v \del (s,\xi_s,\mu_s)(\tilde{X_s})\cdot \tilde{a_s}\right) \, ds \\ \notag &\quad+ \int_0^t \partial_x u(s,\xi_s,\mu_s)\cdot(\gamma_s \, dB_s),
	\end{align} 
	where $(\tilde{\Omega},\tilde{\FF}, \tilde{\P})$ is a copy of $(\Omega, \FF,\P)$ and $(\tilde{X},\tilde{b}, \tilde{\sigma})$ is an independent copy of $(X,b,\sigma)$. 
\end{Thm}

Let us now give examples of functions belonging to the space $\mathcal{W}_2(\R^d).$

\begin{Ex}\label{exlinear2}
	Let $g \in \CC^0(\R^{2d};\R)$ be a function such that its distributional derivative $\nabla g$ belongs to $(W^{1,k}(\R^{2d}))^{2d}$ for some $k \geq 5d.$  Then, the function $$u: \left\{ \begin{array}{rll}
	\R^d \times \PPP_2(\R^d) &\rightarrow \R \\ (x,\mu)&\mapsto \displaystyle\int_{\R^d} g(x,y) \, d\mu(y) \end{array} \right.$$ belongs to $\mathcal{W}_2(\R^d)$ for $d_{\PP}=d_k.$
\end{Ex}

The proof is postponed to the Appendix (Section \ref{proofexlinear2}).

\begin{Ex}\label{ex1}
	Let $F \in \CC^1(\R^d\times \R ; \R)$ be a function such that for all $R>0$  $$ y \in \R \mapsto \nabla F(\cdot,y) \in (W^{1,k_1}(B_R))^{d+1},$$ is well-defined and continuous for some $k_1 \geq d+1.$ Let $g \in \CC^0(\R^d;\R)$ be such that the distributional derivative $\nabla g$ belongs to $(W^{1,k_2}(\R^d))^d$ for some $k_2 \geq 2d.$ Then $$u: \left\{ \begin{array}{cll}
	\R^d \times \PPP_2(\R^d) &\rightarrow \R \\ (x,\mu)&\mapsto F\left(x,\int_{\R^d} g \, d\mu \right)  \end{array} \right.$$ belongs to $\mathcal{W}_2(\R^d)$ for $d_{\PP}=W_2.$  
\end{Ex}
The proof is again postponed to the Appendix (Section \ref{proofex1}).

\begin{Rq}In the abstract, we said that our Itô-Krylov's formula for a flow of measure was the almost analogue of the standard Itô-Krylov formula. We used the word "almost" because Assumption $(1)$ in Definition \ref{sobolevW} is not completely satisfactory. Indeed, we do not assume Sobolev regularity with respect to time, as it is the case in Itô-Krylov's formula for functions defined on $[0,T]\times \R^d.$  Of course if $u$ is of the form $u(t,\mu)=\int_{\R^d} g(t,x) \, d\mu(x)$ with $g\in \CC^0([0,T]\times \R^d;\R)$ at most of quadratic growth in $x$ uniformly in $t$, and such that the distributional derivatives $\partial_t g,$ $\partial_x g$ and $ \partial_x^2 g$ are in $L^k([0,T]\times \R^d)$ for some $k \geq d+1,$ we will succeed in proving Itô-Krylov's formula for $u.$\\
	
	 Let us give the idea of the proof. We regularize $u$ by setting $u^n(t,\mu) := \int_{\R^d} g*\rho_n(t,x) \, d\mu(x),$ where $(\rho_n)_n$ is a mollifying sequence on $\R \times \R^d.$ The function $u^n$ clearly satisfies the assumptions of the standard Itô formula for a flow of measures (see Proposition 5.102 in \cite{CarmonaProbabilisticTheoryMean2018}). It ensures that  for all $t \in [0,T]$ \begin{align}\label{regularityintime}
		\notag u^n(t,\mu_t) &= u^n(0,\mu_0) + \int_0^t \E (\partial_t g * \rho_n (s,X_s)) \, ds + \int_0^t \E \left( \partial_x g * \rho_n(s,X_s) \cdot b_s \right) \, ds \\ &\quad+ \frac{1}{2} \int_0^t \E \left( \partial^2_x g * \rho_n(s,X_s) \cdot a_s \right) \, ds.
		\end{align}  As $g$ is continuous, $(g * \rho_n)_n$ converges to $g$ uniformly on compact sets. It follows from the growth assumption on $g$ that $u^n$ converges point-wise to $u$. Using that $(\partial_t g* \rho_n)_n$ converges in $L^k([0,T]\times \R^d)$ to $\partial_t g$ as $n \rightarrow + \infty,$ we deduce with Krylov's inequality in Corollary \ref{corkrylov} that for all $t \in [0,T]$ $$ \int_0^t \E (\partial_t g * \rho_n (s,X_s)) \, ds  \rightarrow \int_0^t \E (\partial_t g (s,X_s)) \, ds.$$ The same holds with the two other integrals in \eqref{regularityintime}. Taking the limit $n \rightarrow + \infty$ in \eqref{regularityintime} yields for all $t \in [0,T]$ \begin{align*}
		 u(t,\mu_t) &= u(0,\mu_0) + \int_0^t \E (\partial_t g (s,X_s)) \, ds + \int_0^t \E \left( \partial_x g(s,X_s) \cdot b_s \right) \, ds \\ &\quad+ \frac{1}{2} \int_0^t \E \left( \partial^2_x g (s,X_s) \cdot a_s \right) \, ds.
		\end{align*} 
	In the general case, when the dependence in $\mu$ of the function $u$ is not explicit, we cannot apply Krylov's inequality. Indeed, consider a function $ u: [0,T] \times \PPP_2(\R^d) \rightarrow \R$ such that, for all $\mu \in \PPP_2(\R^d),$ $u(\cdot,\mu) \in W^{1,k}([0,T]).$ In Itô's formula for $u$, as in the classical formula, there should be the term $ \int_0^t \partial_t u(s,\mu_s) \, ds$. The assumption does not imply that this term is well-defined. One possible hypothesis is to assume that for all compact $\KK \subset \PPP_2(\R^d),$ $\sup_{\mu \in \KK} |\partial_tu(\cdot,\mu)| \in L^1([0,T]).$ Following our strategy to prove Itô-Krylov's formula, we would consider the mollified version of $u$ defined by $ u^n(t,\mu) := u(\cdot,\mu* \rho^1_n)*\rho^2_n(t),$ where $(\rho^1_n)_n$ and $(\rho^2_n)_n$ are mollifying sequences on $\R^d$ and on $\R$  respectively. Assume that we have proved Itô's formula for $u^n.$ In order to take the limit and deduce Itô's formula for $u,$  we would like to show that \begin{equation*}\label{hypothesistime}
	 \int_0^T |\partial_t u(\cdot,\mu_s*\rho_n^1)*\rho_n^2(s) - \partial_t u(s,\mu_s)| \, ds \rightarrow0. \end{equation*}  However, this convergence is not obvious in the general case since the presence of $\mu_s$ prevents us from using the classical results on convolution and we cannot apply Krylov's inequality if the dependence in the measure argument is not linear.  
	 
	\end{Rq}

	\section{Preliminaries}\label{sectionpreliminaries}

	\subsection{Krylov's inequality and densities.}

	The key element to prove the theorem is Krylov's inequality. We recall it in the next theorem taken from \cite{KrylovControlledDiffusionProcesses2009} (see Theorem $4$ in Section $2.3$). 
	
	\begin{Thm}[Krylov's inequality]\label{inegalitekrylov}
		Let $b: \R^+ \times \Omega \rightarrow \R^d$ and $ \sigma: \R^+ \times \Omega \rightarrow \R^{d \times d_1}$ be two progressively measurable functions. We assume that $p,d_1 \geq d$. Moreover, assume that there exists $K>0$ and $\delta >0$ such that \begin{enumerate}
			\item[] \textbf{(A1)} $\forall (t,\omega) \in \R^+ \times \Omega, \, |b_t(\omega)| + | \sigma_t (\omega)| \leq K$
			\item[] \textbf{(A2)} $\forall (t,\omega) \in \R^+ \times \Omega, \, \forall \lambda \in \R^d, \,  a_t(\omega) \lambda \cdot \lambda  \geq \delta  |\lambda|^2,$ where $a = \sigma \sigma^*$. 
		\end{enumerate}
		For $X_0$ a $\R^d$-valued $\FF_0$-measurable random variable, we define the Itô process $X=(X_t)_t,$ for all $t \in [0,T],$ by  $$  X_t = X_0 + \int_0^t b_s \, ds + \int_0^t \sigma_s \, dB_s.$$
		Let $\lambda >0$ be a positive constant. Then, there exists a constant $N = N( d,p,\lambda, \delta, K)$ such that for all measurable function $f: \R^+ \times \R^d \rightarrow \R$  $$\E \displaystyle\int_0^{\infty} e^{-\lambda t} |f(t,X_t)| \, dt \leq N \Vert f \Vert_{L^{p+1}(\R^+ \times \R^d)}.$$
	\end{Thm}
	
 We will use the following corollary for a finite horizon of time.

	\begin{Cor}\label{corkrylov}
		If $b$ and $\sigma$ satisfy Assumptions \textbf{(A)} and \textbf{(B)}, there exists $N_1 = N_1(d,p,\delta, K,T)$ such that for all measurable function $f: [0,T] \times \R^d \rightarrow \R,$ we have  $$ \E \int_0^T |f(s,X_s)| \, ds \leq  N_1 \Vert f \Vert_{L^{p+1}([0,T] \times \R^d)}.$$ 
	\end{Cor}
	
	\begin{dem}	 We set $b_t=b_T$ and $\sigma_t = \sigma_T$ for $t > T$ to guarantee that Assumptions \textbf{(A1)} and \textbf{(A2)} are satisfied, without changing the process $X$ on $[0,T].$ It remains to apply Krylov's inequality to $\tilde{f}(t,x) := f(t,x) \1_{ t\in[0,T]},$ which gives the existence of $N_1=N_1(d,p,\delta,K)$ such that   $$ e^{-T} \E \int_0^T |f(s,X_s)| \, ds \leq  N_1 \Vert f \Vert_{L^{p+1}([0,T] \times \R^d)}.$$  \end{dem}
	
Krylov's inequality also provides the existence of a density with respect to the Lebesgue measure for $\mu_s$, for almost all $s \in [0,T]$. 

\begin{Prop}\label{densiteexistence}
	Under Assumptions $\textbf{(A)}$ and $\textbf{(B)}$ on the coefficients $b$ and $\sigma$, there exists a function $p \in L^{1}([0,T] \times \R^d;\R^+)\cap L^{(d+1)'}([0,T]\times \R^d;\R^+)$ such that for all $f:[0,T]\times \R^d \rightarrow \R^+$ measurable \begin{equation}
		 \int_0^T \E f(s,X_s) \, ds = \int_{[0,T] \times \R^d} f(s,x) p(s,x) \, dx\,ds.
		 \end{equation}
	If $\tau$ is a stopping time such that $(X_t)_{t \in [0,T]}$ belongs to $B_R$ almost surely on the set $\{ \tau >0\}$, then \begin{equation}\label{inegkrylovdensity}
		\E \int_0^{\tau \wedge T}  f(s,X_s) \, ds \leq \int_{[0,T] \times B_R} f(s,x) p(s,x) \, dx\,ds.\end{equation} Moreover, for almost all $s \in [0,T]$, $\mu_s = \LL(X_s)$ is equal to $p(s,\cdot) \,dx.$
\end{Prop}

We give the proof for the sake of completeness.

\begin{dem}
 We denote by $\mu$ the push-forward measure of $\lambda \otimes \P$, where $\lambda$ is the Lebesgue measure on $[0,T],$ by the measurable map $(t,\omega) \in [0,T] \times \Omega \mapsto (t,X_t(\omega)) \in [0,T] \times \R^d$ defined, for any $A \in \BB([0,T])\otimes\BB(\R^d),$ by  $$ \mu(A) = \int_0^T \E \1_A(s,X_s) \, ds.$$ Note that $\mu$ is a finite measure on $[0,T]\times \R^d.$ The monotone convergence theorem and Krylov's inequality ensure that for all $f:[0,T]\times \R^d \rightarrow \R^+$ measurable  $$ \int_0^T \E f(s,X_s) \, ds = \int_{[0,T]\times \R^d} f(s,x) \, d\mu(s,x) \leq C\|f\|_{L^{p+1}([0,T] \times \R^d)}.$$  Taking $f= \1_{A},$ for $A \in \BB([0,T])\otimes\BB(\R^d)$ with Lebesgue measure $0$, we deduce that $\mu(A)=0.$ Thus $\mu$ is absolutely continuous with respect to the Lebesgue measure on $[0,T]\times\R^d.$ Radon-Nikodym's theorem provides the existence of $p \in L^1([0,T]\times \R^d;\R^+)$ such that for all measurable function $f:[0,T]\times \R^d \rightarrow \R^+$  \begin{equation}\label{equationdensity}
 	\int_0^T \E f(s,X_s) \, ds = \int_{[0,T] \times \R^d} f(s,x) p(s,x) \, dx\,ds. \end{equation} Krylov's inequality exactly proves that the map  $ f \in L^{d+1}([0,T]\times \R^d) \mapsto \displaystyle\int_{[0,T] \times \R^d} f(s,x) p(s,x) \, dx\,ds$ is a continuous linear form. Since the dual space of $L^{d+1}([0,T]\times \R^d)$ is $L^{(d+1)'}([0,T] \times \R^d),$  $p$ belongs to $ L^{(d+1)'}([0,T]\times \R^d).$ \\
 
 \noindent To prove \eqref{inegkrylovdensity}, it is enough to notice that $$ \E \int_0^{\tau \wedge T}  f(s,X_s) \, ds \leq \E \int_0^T f(s,X_s) \1_{B_R}(X_s) \, ds.$$
 Next, we establish that for almost all $s \in [0,T]$, $\mu_s = p(s,\cdot) \,dx.$ We fix $s \in [0,T]$, $n\geq 1$ large enough and $A \in \BB(\R^d).$ Applying \eqref{equationdensity} with $f= \1_{[s-1/n,s+1/n]\times A},$ and using Fubini-Tonelli's theorem, we deduce that $$ \frac n2\int_{s-1/n}^{s+1/n} \P(X_t\in A) \, dt = \frac n2 \int_{s-1/n}^{s+1/n}\int_A p(t,x) \, dx \, ds.$$ Since $t \mapsto \P(X_t \in A)$ is bounded and as Fubini's theorem implies that $ t \mapsto \int_A p(t,x)\, dx$ belongs to $L^1([0,T]),$ it follows from Lebesgue differentiation theorem (see Theorem 7.7 in \cite{RudinRealcomplexanalysis1987}) that for almost all $s \in [0,T]$  $$ \P(X_s \in A) = \int_A p(s,x) \, dx.$$ We denote by $\mathcal{R}$ the set of all Borel sets in $\R^d$ of the form $ \prod_{i=1}^{d}]a_i,b_i[,$ with $a_i < b_i$ two rational numbers for all $i$. The set $\mathcal{R}$ is at most countable, thus for almost $s \in [0,T]$  $$ \forall A \in \mathcal{R}, \quad \P(X_s \in A) = \int_A p(s,x) \, dx.$$ The monotone class theorem enables us to conclude.
\end{dem}

Note that for almost all $s \in [0,T]$, $p(s,\cdot) \in L^{(d+1)'}(\R^d)$ using Fubini-Tonelli's theorem. We deduce the following corollary. 

\begin{Cor}\label{cordensite}
	For almost all $s \in [0,T]$, $ \mu_s \in \PP(\R^d).$
\end{Cor}

We now prove two lemmas dealing with the integrability of the density $p.$

\begin{Lemme} \label{integrability1}
	Let $p$ be the density given by Proposition \ref{densiteexistence}. Then for all $k \geq d+1$  $$ s \in [0,T] \mapsto \Vert p(s,\cdot) \Vert_{L^{k'}(\R^d)} \in L^{k/d}([0,T]).$$
\end{Lemme}

\begin{dem}Using Jensen's inequality since  $\frac{k}{k'} = k-1\geq d,$ we obtain that 
	$$\begin{aligned}
		\int_0^T \left(\int_{\R^d} p(s,x)^{k'} \, dx\right)^{\frac{k}{dk'}} \, ds  &=\int_0^T \left(\int_{\R^d} p(s,x)^{k'-1} p(s,x) \, dx\right)^{\frac{k}{dk'}} \, ds \\ &\leq \int_0^T\int_{\R^d} p(s,x)^{\frac{k}{dk'}(k'-1) +1} \, dx \, ds.
	\end{aligned}$$ 
 By definition of the conjugate exponent, we get $$  \int_0^T\int_{\R^d} p(s,x)^{\frac{k}{dk'}(k'-1) +1} \, dx \, ds = \int_0^T\int_{\R^d} p(s,x)^{\frac{1}{d} +1} \, dx \, ds,$$ which is finite since $(d+1)' = \frac{1}{d} +1$ and $p \in L^{(d+1)'}([0,T]\times \R^d).$
\end{dem}

\begin{Lemme}\label{integrability2}
	Let $p$ and $q$ be two densities of two Itô processes of the form \eqref{itoprocess} and satisfying \textbf{(A)} and \textbf{(B)} given by Proposition \ref{densiteexistence}. Then for $k,\alpha \in \N$ such that $ k \geq \max \{d+1, d(\alpha +1)\},$ we have  $$ \int_0^T \Vert p(s,\cdot) \Vert_{L^{k'}(\R^d)}^{\alpha} \Vert q(s,\cdot) \Vert_{L^{k'}(\R^d)} \, ds < + \infty.$$
\end{Lemme}

\begin{dem}
	Owing to Lemma \ref{integrability1}, the function $s \mapsto \Vert q(s,\cdot) \Vert_{L^{k'}(\R^d)}$ belongs to $L^1([0,T]) \cap L^{k/d}([0,T]).$ Using Hölder's inequality, the proof is complete once we prove that $ s\mapsto \Vert p(s,\cdot) \Vert_{L^{k'}(\R^d)}^{\alpha}$ belongs to $L^r([0,T])$ for some $r \geq \left( \frac{k}{d}\right)'.$ Lemma \ref{integrability1} ensures that $ s\mapsto \Vert p(s,\cdot) \Vert_{L^{k'}(\R^d)}^{\alpha} \in L^{\frac{k}{\alpha d}}([0,T])$ thus we have to prove that $ \left( \frac{k}{d}\right)' \leq \frac{k}{\alpha d}.$ This is equivalent to our assumption $k \geq d(\alpha +1).$
\end{dem}

\subsection{Classical results on convolution and regularization.}

	 Fix $p \in [1+\infty[.$ We will need the two following basic lemmas, which we recall for the sake of clarity. 
	
	\begin{Lemme}[Convolution]\label{convolution}

		\begin{enumerate}
			\item[-] For all $f \in L^p(\R^d)$ and for all $g \in L^1(\R^d)$, the convolution $f * g$ is well-defined and belongs to $L^p(\R^d)$. Moreover, we have $\Vert f * g \Vert_{L^p} \leq \Vert f \Vert_{L^p} \Vert g \Vert_{L^1}.$
			
			\item[-]  For all $f \in L^p(\R^d)$ and for all $g \in L^{p'}(\R^d)$, the convolution $f * g$ is well-defined and belongs to $L^{\infty}(\R^d)$. Moreover, we have $\Vert f * g \Vert_{L^{\infty}} \leq \Vert f \Vert_{L^p} \Vert g \Vert_{L^{p'}}.$
		\end{enumerate}
	\end{Lemme}

	\begin{Lemme}[Regularization]\label{regularisation}
		\
		Recall that $(\rho_n)_n$ is a mollifying sequence. 
		\begin{enumerate}
			\item[-] Let $f \in L^1_{\text{loc}}(\R^d)$ and $\rho \in \CC^{\infty}_c(\R^d)$. Then $f * \rho \in \CC^{\infty}(\R^d)$ and $ \forall\alpha \in \N^d,$ $ \partial^{\alpha} ( f * \rho) = f * \partial^{\alpha} \rho.$
			
			\item[-] If $f \in L^p(\R^d)$, then $f * \rho_n \overset{L^p}{\longrightarrow} f$, and if $f \in \CC^0(\R^d)$, $f * \rho_n\rightarrow f$ uniformly on compact sets.
			
			\item[-] If $f \in L^p_{loc}(\R^d)$, then for all $R>0$, $f * \rho_n \rightarrow f $ in $L^p(B_R).$
		\end{enumerate}
	\end{Lemme}
	The following proposition will also be useful. 
	
	\begin{Prop}\label{deriveefaible}
		Let $f \in \CC^0(\R^d)$ be a function admitting distributional derivatives of order $1$ et $2$ in $L^1_{loc}(\R^d).$ Then $f * \rho_n \in \CC^{\infty}(\R^d)$ and for all $i,j \in \{1, \dots d \}$  $$ \left\{ \begin{array}{rll} \partial_{x_i} ( f * \rho_n) &= \partial_{x_i}f * \rho_n \\ \partial_{x_i \, x_j} ( f * \rho_n) &= \partial_{x_i \, x_j}f * \rho_n.
		\end{array}\right.$$
	\end{Prop}

The next lemma deals with the convolution of a function $f \in L^p$  with $\mu \in \PPP(\R^d).$
	
\begin{Lemme}\label{continuiteconvolution}
	Let $f \in L^{p}(\R^d)$. Then $\mu \in \PPP(\R^d) \mapsto f * \mu \in L^p(\R^d)$ is continuous.
\end{Lemme}

\begin{dem}Note that the convolution $f * \mu$ is well-defined as an element of $L^p(\R^d)$ thanks to Jensen's inequality which shows that $$ \forall f \in L^p(\R^d), \, \forall \mu \in \PPP(\R^d),\, \Vert f * \mu \Vert_{L^p} \leq \Vert f \Vert_{L^p}.$$
	Let $(\mu_n)_n$ be a sequence of $\PPP(\R^d)$ weakly convergent to $\mu \in \PPP(\R^d)$. Using Skorokhod's representation theorem (see Theorem $6.7$ in \cite{BillingsleyConvergenceProbabilityMeasures1999}), there exists a probability space $(\Omega',\FF', \P')$, a sequence of random variables $(X_n)_n$ converging $\P'$-almost surely to a random variable $X$ such that, the law of $X_n$ is $\mu_n$ for all $n$ and the law of $X$ if $\mu$. For any $a \in \R^d,$ let us denote by $\tau_a f$ the translation of $f$ defined, for all $x \in \R^d,$ by $\tau_a f(x) := f(x-a).$ Jensen's inequality and Fubini-Tonelli's theorem yield \begin{align*}
	\Vert f * \mu_n - f * \mu \Vert_{L^p}^p &= \int_{\R^d} \left\vert \E'(f(x-X_n) - f(x-X))\right\vert^p \, dx \\ &\leq \int_{\R^d} \E'(\left\vert f(x-X_n) - f(x-X)\right\vert^p) \, dx  \\ &= \E' (\Vert \tau_{X_n-X}f - f \Vert_{L^p}^p).
	\end{align*}
	It follows from the almost sure convergence of $(X_n)_n$ to $X$ and the continuity of the translation operator in $L^p$ that $\Vert \tau_{X_n-X}f - f \Vert_{L^p}^p \overset{a.s.}{\longrightarrow} 0.$ Moreover, the inequality  \begin{align*} \Vert \tau_{X_n-X}f - f \Vert_{L^p}^p &\leq 2^{p-1} ( \Vert \tau_{X_n-X}f \Vert_{L^p}^p + \Vert f \Vert_{L^p}^p )\\ &= 2^{p} \Vert f \Vert_{L^p}^p, \end{align*}enables us to conclude with the dominated convergence theorem.
\end{dem}

\subsection{Convolution of probability measures}

\begin{Lemme}[Contraction inequality]\label{convolutionmesure}
	Fix $\mu,\nu,m \in \PPP_2(\R^d)$. Then, we have  $$ W_2(\mu * m, \nu * m) \leq W_2 ( \mu, \nu).$$
\end{Lemme}

\begin{dem}Let $\pi \in \PPP_2(\R^d \times \R^d)$ be an optimal coupling between $\mu$ and $\nu$. 
	We consider a couple of random variables $(X,Y)$ with law $\pi,$ and a random variable $Z$ independent of $(X,Y)$ with law $m$. The law of $X+Z$ being $\mu * m$ and the law of $Y+Z$ being $\nu * m$, one has \begin{align*}
	W_2(\mu * m, \nu * m) \leq \Vert (X+Z)-(Y+Z) \Vert_{L^2}  = W_2(\mu,\nu).
	\end{align*}
\end{dem}

The next corollary follows from the fact that $\rho_n \overset{W_2}{\longrightarrow} \delta_0.$

\begin{Cor}\label{regularisationconvolutionmeasure}
	For all $\mu \in \PPP_2(\R^d),$ $ \mu * \rho_n \overset{W_2}{\longrightarrow} \mu.$
\end{Cor}

	\subsection{Measurability}

	We will need the following lemma to guarantee that, for $u \in \mathcal{W}_1(\R^d),$ we can find versions of $\partial_v \del$ and $\partial^2_v \del$ that are measurable with respect to $(\mu,v) \in \PP(\R^d) \times \R^d$. 
	
	\begin{Lemme}\label{measurabilty}
		Let $u: E \rightarrow L^{k}(\R^d)$ be a continuous function, where $E$ is a metric space and $k >1$. Then, for all $x \in E$, we can find a version of $u(x)$ such that $(x,v) \in E \times \R^d \mapsto u(x)(v)$ is measurable with respect to $\BB(E)\otimes \BB(\R^d).$ 
	\end{Lemme}
	
	\begin{dem}
		For $(x,v) \in E \times \R^d$, we define  $$\tilde{u}(x,v) = \underset{n \rightarrow + \infty}{\underline{\lim}} \frac{1}{\lambda(B(v,1/n))} \int_{B(v,1/n)} u(x)(y) \, dy = \underset{n \rightarrow + \infty}{\underline{\lim}} u^n(x,v),$$ where $\lambda$ denotes the Lebesgue measure on $\R^d$. From Lebesgue differentiation theorem (see Theorem 7.7 in \cite{RudinRealcomplexanalysis1987}), we deduce that for all $x \in E$, $\tilde{u}(x,\cdot) = u(x)$ $\lambda$-almost everywhere. We prove that for all $n\geq 1$, $u^n$ is continuous. Note that $\frac{1}{\lambda(B(v,1/n))}$ does not depend on $v$. The continuity of $u^n$ follows from the continuity of $x \in E \mapsto u(x) \in L^{k}(\R^d)$, $v \in \R^d \mapsto \1_{B(v,1/n)} \in L^{k'}(\R^d)$ (coming from the dominated convergence theorem), and of $(f,g) \in L^{k}(\R^d) \times L^{k'}(\R^d) \mapsto \int_{\R^d} fg \, dx.$ 
	\end{dem}

	\section{Proof of Theorem \ref{ito_krylov}}\label{sectionproof1}

	The proof will be divided into three parts. Step 1 is dedicated to prove that all the terms in Itô-Krylov's formula \eqref{formulaito1} are well-defined. In Step 2, we regularize $u$ by convolution of the measure argument with a mollifying sequence $(\rho_n)_n.$ The effect of replacing $u(\mu)$ by $u(\mu * \rho_n)$ is that the linear derivative is regularized by convolution, in its space variable. Then, we apply the standard Itô's formula for a flow of measure. We finally take the limit $n \rightarrow + \infty$ in Step $3$ with the help of Krylov's inequality. \\

	\noindent\textbf{Step 1: All the terms in \eqref{formulaito1} are well-defined.}\\
	
	  Let us show that the two integrals in \eqref{formulaito1} are well-defined.\\
	  
	\noindent\textbf{Measurability.} Thanks to Lemma \ref{measurabilty}, we can find a version of $\partial_v \del $ which is measurable with respect to $(\mu,v) \in \PP(\R^d)\times \R^d.$ To conclude, we prove that $s \mapsto \mu_s \in \PP(\R^d)$ is measurable. Indeed if it is the case, the function $(s,\omega) \in [0,T]\times \Omega \mapsto \partial_v \del(\mu_s)(X_s(\omega)).b_s(\omega)$ will be measurable by composition. First, note that $\mu_s \in \PP(\R^d)$ for almost all $s \in [0,T]$ (see Corollary \ref{cordensite}) so we can change $\mu_s$ on a negligible set of times $s$ to ensure that $\mu_s \in \PP(\R^d)$ for all $s \in [0,T].$ But $\mu_s = \displaystyle\lim_{n\rightarrow+\infty} \mu_s* \rho_n$ for $d_{\PP}$ by Assumption \textbf{(H2)} in Definition \ref{defespaceP}. It remains to show that $s \mapsto \mu_s * \rho_n \in \PP(\R^d)$ is continuous and thus mesurable for all $n$. This follows from the continuity of $s \mapsto \mu_s \in \PPP_2(\R^d)$ and also from Assumption \textbf{(H1)} in Definition \ref{defespaceP}. \\
		
\noindent	\textbf{Integrability.} We can omit the coefficients $b$ and $a$ to prove the integrability properties because they are uniformly bounded. Taking advantage from the existence of a density coming from Proposition \ref{densiteexistence}, we have by Hölder's inequality  \begin{align*}
		\int_0^T \E \left| \partial_v \del (\mu_s)(X_s) \right| \, ds & = \int_0^T \int_{\R^d}  \left| \partial_v \del (\mu_s)(x) \right| p(s,x) \, dx\, ds \\ &\leq \int_0^T \left\Vert \partial_v \del(\mu_s)(\cdot)\right\Vert_{L^k(\R^d)} \Vert p(s,\cdot)\Vert_{L^{k'}(\R^d)} \, ds \\ &\leq \int_0^T C \left(1+ \Vert p(s,\cdot) \Vert_{L^{k'}(\R^d)}^{\alpha}\right) \Vert p(s,\cdot)\Vert_{L^{k'}(\R^d)} \, ds, 
		\end{align*}
		for some constant $C$ coming from Assumption $(2)$ in Definition \ref{sobolevW1} because the flow $(\mu_s)_{s\leq T}$ is compact in $\PPP_2(\R^d)$ and belongs to $\PP(\R^d)$ for almost all $s.$  The last bound is finite thanks to Lemma \ref{integrability1} since $ k\geq \max \{d(\alpha+1),d+1\}.$ The same properties hold for the term involving $\partial^2_v \del.$\\

\noindent\textbf{Step 2: Itô's formula for the mollification of $\boldsymbol{u}.$}  \\[1pt]
	
	For $n \geq 1$, we set  $u^n: \mu \in \PPP_2(\R^d) \mapsto u (\mu * \rho_n).$  By standard arguments, for each $ n \geq 1,$ $u^n$ has a linear derivative given by \begin{align*}
	\frac{\delta u^n}{\delta m} (\mu)(v) = \int_{\R^d} \del (\mu * \rho_n)(x) \rho_n(v-x)\, dx  = \del (\mu * \rho_n) * \rho_n (v).
	\end{align*}
	Now, we aim at applying the standard Itô formula for a flow of probability measures (see for example Theorem 5.99 in Chapter $5$ of \cite{CarmonaProbabilisticTheoryMean2018} with the L-derivative) to $u^n$ for a fixed $n \geq 1.$ \\

	\textbf{(i) Regularity of $\boldsymbol{\frac{\delta u^{n}}{\delta m}(\mu)}$ for a fixed $\boldsymbol{\mu \in \PPP_2(\R^d)}.$} Since for all $\mu \in \PPP_2(\R^d)$, $ \mu * \rho_{n} \in \PP(\R^d),$ Proposition \ref{deriveefaible} implies that $\frac{\delta u^{n}}{\delta m}(\mu)(.) \in \CC^{\infty}(\R^d)$ and for all $i,j \in \{1, \dots, d\}$  $$ \partial_{v_i}  \frac{\delta u^{n}}{\delta_m}(\mu)  = \partial_{v_i} \frac{\delta u}{\delta_m}(\mu * \rho_{n}) * \rho_{n} \quad \text{and}\quad \partial_{v_i \, v_j} \frac{\delta u^{n}}{\delta m}(\mu) = \partial_{v_i \, v_j}\del (\mu * \rho_{n}) * \rho_{n}.$$ 
		
	\textbf{(ii) Continuity of $\boldsymbol{\partial_v\frac{\delta u^{n}}{\delta m}}$ and $\boldsymbol{\partial^2_v\frac{\delta u^{n}}{\delta m}}$ with respect to $\boldsymbol{(\mu,v)}.$} Let $ i \in \{1, \dots ,d \}$, $(\mu_m)_m \in \PPP_2(\R^d)^{\N}$ and $(v_m)_m \in (\R^d)^{\N}$ be sequences converging respectively to $\mu$ and $v$. We have  \begin{align*}
		&\left\vert \partial_{v_i} \frac{\delta u^{n}}{\delta m}(\mu_m)(v_m) - \partial_{v_i} \frac{\delta u^{n}}{\delta m}(\mu)(v) \right\vert \\&\leq	\left\vert \partial_{v_i} \frac{\delta u^{n}}{\delta m}(\mu_m)(v_m) - \partial_{v_i} \frac{\delta u^{n}}{\delta m}(\mu)(v_m) \right\vert+ \left\vert \partial_{v_i} \frac{\delta u^{n}}{\delta m}(\mu)(v_m) - \partial_{v_i} \frac{\delta u^{n}}{\delta m}(\mu)(v) \right\vert \\ &=: D_1 + D_2
		\end{align*}
		$D_2$ converges to $0$ when $m \rightarrow + \infty$ by $(i).$ For $D_1$, the convolution inequality $L^k*L^{k'}$ gives that
		\begin{align*}
		D_1 & = \left\vert  \partial_{v_i} \frac{\delta u}{\delta_m}(\mu_m * \rho_{n}) * \rho_{n} (v_m) - \partial_{v_i} \frac{\delta u}{\delta_m}(\mu * \rho_{n}) * \rho_{n}(v_m) \right\vert  \\ &\leq \left\Vert \partial_{v_i} \frac{\delta u}{\delta_m}(\mu_m * \rho_{n})  - \partial_{v_i} \frac{\delta u}{\delta_m}(\mu * \rho_{n}) \right\Vert_{L^{k}}  \Vert \rho_{n} \Vert_{L^{k'}}.
		\end{align*}
		Assumption \textbf{(H1)} in Definition \ref{defespaceP} provides that $ \mu_m * \rho_{n} \overset{d_{\PP}}{\longrightarrow} \mu * \rho_{n}$ when $ m \rightarrow + \infty.$ Finally, using the first assumption in Definition \ref{sobolevW1}, we conclude that $D_1$ converges to $0$ when $m \rightarrow + \infty$. This shows the continuity of $\partial_v\frac{\delta u^{n}}{\delta m}$ on $\PPP_2(\R^d) \times \R^d$. The same reasoning proves the joint continuity of $\partial^2_v\frac{\delta u^{n}}{\delta m}$.\\ 
		
	 \textbf{(iii) Boundedness of  $\boldsymbol{\partial_v\frac{\delta u^{n}}{\delta m}}$ and $\boldsymbol{\partial^2_v\frac{\delta u^{n}}{\delta m}}.$} Let $\KK \subset \PPP_2(\R^d)$ be a compact set. For $\mu \in \KK$ and $v \in \R^d$, one has
		\begin{align*}
		\left\vert \partial_{v_i}\frac{\delta u^{n}}{\delta m}(\mu)(v) \right\vert &\leq \left\Vert \partial_{v_i} \frac{\delta u}{\delta_m}(\mu * \rho_{n}) \right\Vert_{L^k}  \Vert \rho_{n} \Vert_{L^{k'}}.
		\end{align*}
		The set $ \{\mu * \rho_{n}, \, \mu \in \KK \}$ is compact in $(\PP(\R^d),d_{\PP})$ as the image of the compact $\KK$ by the application $ \mu \in \PPP_2(\R^d) \mapsto \mu * \rho_{n} \in \PP(\R^d)$ which is continuous by Assumption (\textbf{H1}) in Definition \ref{defespaceP}. The first assumption in Definition \ref{sobolevW1} guarantees that  $ \sup_{\mu \in \KK} \left\Vert \partial_{v_i} \frac{\delta u}{\delta_m}(\mu * \rho_{n})\right\Vert_{L^{k}(\R^d)} < + \infty $ and thus  $$ \underset{v \in \R^d}{\sup}\sup_{\mu \in \KK}  \left\vert \partial_{v}\frac{\delta u^{n}}{\delta m}(\mu)(v) \right\vert< \infty.$$ The same property holds for $\partial^2_{v}\frac{\delta u^{n}}{\delta m}.$\\

	\noindent We can thus apply Itô's formula of \cite{CarmonaProbabilisticTheoryMean2018} to obtain that for all $n \geq 1$ and for all $t \in [0,T]$  \begin{equation}\label{formulau_nito1}
		 u^n(\mu_t) = u^n(\mu_0) + \int_0^t \E \left( \partial_v \frac{\delta u^n}{\delta m} (\mu_s)(X_s)\cdot b_s\right) \,ds + \frac{1}{2}  \int_0^t \E \left( \partial^2_v \frac{\delta u^n}{\delta m}(\mu_s)(X_s)\cdot a_s\right) \,ds.\end{equation}

	\noindent \textbf{Step 3: Letting $\boldsymbol{n \rightarrow + \infty}$.}  \\[1pt]
Our aim is now to take the limit $n \rightarrow + \infty$ in \eqref{formulau_nito1}.	
As for all $\mu \in \PPP_2(\R^d),$ $\mu * \rho_n \overset{W_2}{\longrightarrow} \mu$ and $u$ is continuous on $\PPP_2(\R^d),$ we deduce that $(u^n)_n$ converges pointwise to $u$. It remains to take the limit in the two integrals of \eqref{formulau_nito1}. We show that  \begin{equation}\label{term1Ito1}
		  \int_0^t \E \left( \partial_v \frac{\delta u^n}{\delta m} (\mu_s)(X_s)\cdot b_s\right) \,ds \rightarrow \int_0^t \E \left( \partial_v \frac{\delta u}{\delta m} (\mu_s)(X_s)\cdot b_s\right) \,ds.\end{equation} Since $b$ is uniformly bounded, it is enough to prove that $$ \E \int_0^T \left| \partial_v \del (\mu * \rho_n)* \rho_n (X_s) - \partial_v \del (\mu_s)(X_s) \right| \, ds \rightarrow 0.$$ By Proposition \ref{densiteexistence}, Hölder's inequality and then the $L^1*L^k$ convolution inequality, one has \begin{align*}
		&\E \int_0^T \left| \partial_v\del (\mu_s * \rho_n)* \rho_n (X_s) - \partial_v \del (\mu_s)(X_s) \right| \, ds   \\ &\leq   \int_0^T \left\Vert   \partial_v \del (\mu_s * \rho_n)  - \partial_v \del (\mu_s) \right\Vert_{L^k(\R^d)} \Vert p(s,\cdot) \Vert_{L^{k'}(\R^d)} \, ds  \\ &\quad+  \int_0^T \left\Vert   \partial_v \del (\mu_s)* \rho_n  - \partial_v \del (\mu_s) \right\Vert_{L^k(\R^d)} \Vert p(s,\cdot) \Vert_{L^{k'}(\R^d)} \, ds\\ &=: I_1 + I_2.
	\end{align*}
 The integrand in $I_1$ converges to $0$ for almost all $s$ using Assumption (1) in Theorem \ref{ito_krylov} and the fact that $ \mu_s * \rho_n \overset{d_{\PP}}{\longrightarrow} \mu_s$ for almost all $s$ thanks to Assumption \textbf{(H2)} in Definition \ref{defespaceP}. Let us now prove that the dominated convergence theorem applies. The integrand is bounded by  $$ \left[\sup_{n\geq 1} \left\Vert   \partial_v \del (\mu_s * \rho_n) \right\Vert_{L^k(\R^d)} +   \left\Vert   \partial_v \del (\mu_s ) \right\Vert_{L^k(\R^d)}\right] \Vert p(s,\cdot) \Vert_{L^{k'}(\R^d)}.$$ 	Note that the set $\{\mu_s * \rho_n, \, s \in [0,T], \, n \geq 1 \} \cup \{ \mu_s, \, s \in [0,T] \}$ is compact in $\PPP_2(\R^d).$ Indeed, if $(s_k)_k \in [0,T]^{\N}$ and $(n_k)_k \in \N^{\N}$ are two sequences, we have to find a convergent subsequence from $(\mu_{s_k}* \rho_{n_k})_k$. Up to an extraction, we can assume that  $(s_k)_k$ converges to some $s \in [0,T]$. There are two cases. If there exists $l$ such that $n_k =l$ infinitely often, then $\mu_{s_k} * \rho_l \overset{W_2}{\longrightarrow} \mu_s * \rho_l$ by the contraction inequality (see Lemma \ref{convolutionmesure}).
 Otherwise, we can assume that $(n_k)_k$ converges to $+ \infty$. We use the triangle inequality to get \begin{align*}
 W_2 (\mu_{s_k} * \rho_{n_k}, \mu_s) &\leq W_2(\mu_{s_k} * \rho_{n_k},\mu_{s} * \rho_{n_k}) + W_2(\mu_{s} * \rho_{n_k}, \mu_s).
 \end{align*}
 The last term converges to $0$ owing to Lemma \ref{regularisationconvolutionmeasure}, and the first is bounded by $W_2(\mu_{s_k},\mu_s)$ by the contraction inequality (see Lemma \ref{convolutionmesure}), which converges to $0.$ Thus Assumption $(2)$ in Definition \ref{sobolevW1} ensures that there exists $C>0$ such that for almost all $s \in [0,T]$ and for all $n$  $$\left\Vert   \partial_v \del (\mu_s * \rho_n) \right\Vert_{L^k(\R^d)} +   \left\Vert   \partial_v \del (\mu_s ) \right\Vert_{L^k(\R^d)} \leq C (1+\Vert p(s,\cdot)*\rho_n \Vert_{L^{k'}(\R^d)}^{\alpha} + \Vert p(s,\cdot) \Vert_{L^{k'}(\R^d)}^{\alpha}).$$ It follows from the convolution inequality $L^{k'}*L^1$ that for almost all $s$  \begin{align*}
 	  &\left[\sup_{n\geq 1} \left\Vert   \partial_v \del (\mu_s * \rho_n) \right\Vert_{L^k(\R^d)} +   \left\Vert   \partial_v \del (\mu_s ) \right\Vert_{L^k(\R^d)}\right] \Vert p(s,\cdot) \Vert_{L^{k'}(\R^d)} \\ &\leq 2C (1 + \Vert p(s,\cdot) \Vert_{L^{k'}(\R^d)}^{\alpha} )\Vert p(s,\cdot) \Vert_{L^{k'}(\R^d)}, \end{align*} which is integrable on $[0,T]$ thanks to Lemma \ref{integrability1} since $ k\geq \max \{d(\alpha+1),d+1\}.$
 We conclude by the dominated convergence theorem that $I_1$ converges to $0$. The term $I_2$ also converges to $0$ following the same method. Indeed, for almost all $s$, $\partial_v\del(\mu_s)(\cdot) \in L^k(\R^d)$ thus the integrand converges to $0$ by Lemma \ref{regularisation} and we conclude with the dominated convergence theorem. Therefore \eqref{term1Ito1} is proved. Following the same lines, we take the limit $n \rightarrow + \infty$ in the last integral of \eqref{formulau_nito1} to obtain that for all $t \in [0,T]$   $$ \int_0^t \E \left( \partial^2_v \frac{\delta u^n}{\delta m} (\mu_s)(X_s)\cdot a_s\right) \,ds \rightarrow \int_0^t \E \left( \partial^2_v \frac{\delta u}{\delta m} (\mu_s)(X_s)\cdot a_s\right) \,ds.$$ This concludes the proof of Theorem \ref{ito_krylov}. \hfill$\square$
	
	\section{Proof of Theorem \ref{extensionito}}\label{sectionproof2}
	
		The strategy of the proof is the following. In Step 1, we prove some integrability results coming from Assumption $(5)$ in Definition \ref{sobolevW}. Step 2 is devoted to prove that all the terms in Itô-Krylov's formula \eqref{formulaito2} are well-defined using a localization argument, Krylov's inequality, and Step 1. Moreover, we see that it is enough to prove the formula up to random times localizing the process $\xi.$ Step 3 is dedicated to regularize $u$ using convolutions both in space and measure variables. In Step 4 and 5, we follow the strategy of the proof of Theorem $5.102$ in \cite{CarmonaProbabilisticTheoryMean2018} to prove Itô-Krylov's formula for $u^n$, the mollified version of $u.$ Finally, Step 6 aims at taking the limit $ n \rightarrow + \infty$ thanks to Krylov's inequality.\\
		
		 Note that there are three kind of integrals in Itô's formula \eqref{formulaito2}: the terms involving standard time and space derivatives in the first line, those involving the linear derivative in the second line and the martingale term in the third line. We will treat them separately.   \\
		
		\noindent\textbf{Step 1: Useful integrability results.}\\
		
	It follows from Assumption $(5)$ in Definition \ref{sobolevW} and Lemma \ref{integrability2} that for any $M>0$ the following quantities are finite:
		
		\begin{align}
			\label{bound1}&J_1(M):= \displaystyle\int_0^T \left[\sup_{n \geq 1} \left\Vert \partial_x u(s,\cdot,\mu_s* \rho_n) \right\Vert_{L^{k_1}(B_M)} + \sup_{n \geq 1}\left\Vert\partial^2_x u  (s,\cdot,\mu_s* \rho_n) \right\Vert_{L^{k_1}(B_M)} \right] \Vert q(s,\cdot)\Vert_{L^{k_1'}(B_M)} \, ds,\\ \notag&J_2(M):=\displaystyle\int_0^T \sup_{n \geq 1} \left\Vert \partial_x u(s,\cdot,\mu_s* \rho_n) \right\Vert_{L^{2k_1}(B_M)}^2  \Vert q(s,\cdot)\Vert_{L^{k_1'}(B_M)} \, ds,\\\notag&J_3(M):=\displaystyle\int_0^T \sup_{n \geq 1} \left\Vert \partial_v \del(s,\cdot,\mu_s* \rho_n)(\cdot) \right\Vert_{L^{k_2}(B_M\times\R^{d})} \Vert q(s,\cdot)\Vert_{L^{k_2'}(B_M)}  \Vert p(s,\cdot)\Vert_{L^{k_2'}(\R^d)} \, ds,\\\notag&J_4(M):=\displaystyle\int_0^T  \sup_{n \geq 1}\left\Vert\partial^2_v \del  (s,\cdot,\mu_s* \rho_n)(\cdot) \right\Vert_{L^{k_2}(B_M\times\R^{d})}  \Vert q(s,\cdot)\Vert_{L^{k_2'}(B_M)}  \Vert p(s,\cdot)\Vert_{L^{k_2'}(\R^d)} \, ds.
		\end{align}

	To prove this, we follow the method employed in Step $3$ of the preceding proof to justify the dominated convergence theorem. We just give details for $J_2(M)$ since it requires a bit more attention. Owing to Assumption $(2)$ in Definition \ref{sobolevW}, we know that for all $(t,\mu) \in [0,T] \times \PP(\R^d),$ $\partial_x u(t,\cdot,\mu) \in W^{1,k_1}(B).$ Sobolev embedding theorem (see Corollary 9.14 in \cite{BrezisFunctionalAnalysisSobolev2010}) ensures that the embedding $  W^{1,k_1}(B_M) \hookrightarrow L^{2k_1}(B_M)$ is continuous since $ k_1 \geq d+1.$ Thus there exists $C>0$ such that  $$ \forall t \in [0,T], \, \forall \mu \in \PP(\R^d),\, \Vert \partial_x u(t,\cdot,\mu) \Vert_{L^{2k_1}(B_M)} \leq C \left(\Vert\partial_x u(t,\cdot,\mu) \Vert_{L^{k_1}(B_M)} + \Vert\partial_x^2 u(t,\cdot,\mu) \Vert_{L^{k_1}(B_M)}\right).$$ Thanks to Assumption $(5)$ in Definition \ref{sobolevW}, there exists a constant $C_{M}>0$ such that for almost all $s$ and for all $n\geq 1$  $$ \sup_{n\geq 1} \Vert \partial_x u(s,\cdot,\mu_s* \rho_n) \Vert^2_{L^{2k_1}(B_M)} \leq C_{M} \left(1 + \Vert p(s,\cdot) \Vert_{L^{k_1'}(\R^d)}^{2\alpha_1}\right),$$ where we used the fact that $\{ \mu_s * \rho_n, \, s \in [0,T], \, n \geq 1 \}$ is relatively compact in $\PPP_2(\R^d)$ and the convolution inequality $L^{k_1'}*L^1.$ We conclude with Lemma \ref{integrability2} since $k_1 \geq \max \{ d(2\alpha_1+1),d+1\}.$ Note that these integrability properties remain true if we replace $\mu_s * \rho_n$ by $\mu_s$ and remove the supremum. We justify it only for the second point. It follows from the continuity assumption $(2)$ in Definition \ref{sobolevW} that  for almost all $s \in [0,T]$  $$  \partial_x u(s,\cdot,\mu_s* \rho_n)  \overset{W^{1,k_1}(B_M)}{\longrightarrow} \partial_x u(s,\cdot,\mu_s),$$ because  $\mu_s * \rho_n \overset{d_{\PP}}{\longrightarrow} \mu_s$ for almost all $s$. Sobolev embedding theorem guarantees that  $$ \Vert \partial_x u(t,\cdot,\mu* \rho_n) \Vert_{L^{2k_1}(B_M)} \rightarrow  \Vert \partial_x u(t,\cdot,\mu) \Vert_{L^{2k_1}(B_M)}.$$ Thus we obtain  \begin{equation*}
		 \int_0^T \Vert \partial_x u(s,\cdot,\mu_s) \Vert_{L^{2k_1}(B_M)}^2 \Vert q(s,\cdot)\Vert_{L^{k_1'}(B_M)} \, ds   \leq J_2(M) <+ \infty.\end{equation*} \\

	\noindent\textbf{Step 2: Meaning of the terms in \eqref{formulaito2} and localization.}\\
	
 Let $(T_M)_M$ be the sequence of stopping times converging almost surely to $T$ defined by  $$T_M= \inf \{ t \in[0,T],\, |\xi_t| \geq M \}\wedge T.$$ Let $\xi^M_t = \xi_{t\wedge T_M},$ which is bounded by $M$ on the set $\{T_M > 0 \}$. \\
 
 \noindent\textbf{(i) Terms involving standard derivatives in \eqref{formulaito2}.} We prove that almost surely $$\int_0^T |\partial_x u(s,\xi_s,\mu_s)\cdot\eta_s| \, ds < + \infty.$$ By Proposition \ref{densiteexistence} and Hölder's inequality, one has  \begin{align*}\E \int_0^{T\wedge T_M} |\partial_x u (s, \xi_s,\mu_s)| \,ds &\leq  \int_0^T \int_{B_{M}}  |\partial_x u (s, x,\mu_s)|q(s,x) \, dx\, ds\\ &\leq \int_0^T \Vert  \partial_x u (s, \cdot,\mu_s)\Vert_{L^{k_1}(B_M)} \Vert q(s,\cdot) \Vert_{L^{k_1'}(B_M)} \, ds \\ &\leq J_1(M),
		\end{align*}which is finite (see \eqref{bound1} in Step $1$). We deduce that almost surely, for all $ M \geq 1$  $$\int_0^{T\wedge T_M} |\partial_x u (s, \xi_s,\mu_s)| \,ds < \infty.$$ But it is clear that for almost all $\omega \in \Omega$ and for $M$ bigger than some random constant $M(\omega) \geq 1$, $T_M(\omega) =T.$ Thus, since $\eta$ is uniformly bounded, $\int_0^T |\partial_x u(s,\xi_s,\mu_s).\eta_s| \, ds$ is finite almost surely. The other terms in the first line of Itô's formula \eqref{formulaito2} are treated with the same method.\\
		
	\noindent \textbf{(ii) Martingale term in \eqref{formulaito2}.}  We need to prove that $\int_0^T |\partial_x u(s,\xi_s,\mu_s) |^2 \, ds$ is almost surely finite. Reasoning as before, it is a consequence of the fact that $J_2$ is finite since we have $$  \int_0^T \Vert  \partial_x u (s, \cdot,\mu_s)\Vert_{L^{2k_1}(B_M)}^2 \Vert q(s,\cdot) \Vert_{L^{k_1'}(B_M)} \, ds \leq J_2(M). $$  Therefore the martingale term in \eqref{formulaito2} is well-defined.\\
	
	\noindent \textbf{(iii) Terms involving the linear derivative in \eqref{formulaito2}.} We remark that $\tilde{X}$ and $\xi$ can be seen as independent processes on the product space $\Omega \times \tilde{\Omega}$ with $\LL(\tilde{X}_s) = p(s,\cdot) \, dx$ and $\LL(\xi_s) = q(s,\cdot) \, dx$  for almost all $s.$ Hölder's inequality gives that \begin{align*}&\E\int_0^{T \wedge T_M} \tilde{\E} \left|\partial_v \del (s,\xi_s,\mu_s)(\tilde{X_s})\right| \, ds \\&\leq  \int_0^T \int_{B_M \times \R^d} \left|\partial_v \del (s,x,\mu_s)(v)\right|q(s,x) p(s,v) \, dx \, dv \, ds \\ &\leq \int_0^T \left\Vert \partial_v \del (s,\cdot,\mu_s)(\cdot)\right\Vert_{L^{k_2}(B_M\times \R^d)} \Vert q(s,\cdot) \Vert_{L^{k_2'}(B_M)}  \Vert p(s,\cdot) \Vert_{L^{k_2'}(\R^d)} \, ds \\ &=J_3(M),\end{align*} 
		  which was defined in \eqref{bound1} and is finite. We deduce as previously that $\int_0^T \tilde{\E} \left|\partial_v \del (s,\xi_s,\mu_s)(\tilde{X_s}).\tilde{b_s}\right| \, ds$ is almost surely finite. The term involving $\partial_v^2 \del$ is dealt similarly. \\
	
	 Since all the terms in \eqref{formulaito2} are well-defined, it is enough to prove Itô-Krylov's formula for 
	 \newline $u(t\wedge T_M, \xi_{t \wedge T_M}, \mu_{t\wedge T_M})$ almost surely for all $t \in [0,T],$ and then take the limit $ M \rightarrow + \infty$ using the continuity of the integrals in Itô-Krylov's formula with respect to $t$. So we fix $\tau := T_M$ for $M \geq 1$ and we want to prove the formula up to time $\tau.$\\
	
	\noindent\textbf{Step 3: Mollification of $\boldsymbol{u}$.}\\
	
	Let $u^n$ be the function defined by $u^n(t,x,\mu):= u(t,\cdot,\mu* \rho_n)*\rho_n(x).$ It is clearly continuous on $[0,T]\times \R^d \times \PPP_2(\R^d),$ as $u.$ Since $\partial_t u$ is jointly continuous, it follows from Leibniz's rule that $u^n$ is $\CC^1$ with respect to $t$ and that we can differentiate under the integral i.e. for all $(t,x,\mu) \in [0,T]\times \R^d \times \PPP_2(\R^d)$  $$ \partial_t u^n(t,x,\mu) = \partial_t u(t,\cdot, \mu * \rho_n) * \rho_n(x),$$ which is also jointly continuous. As a result of Lemma \ref{regularisation} and Proposition \ref{deriveefaible},  $u^n$ is $\CC^2$ with respect to $x$ and we have $$\partial_x u^n(t,x,\mu) = \partial_x u(t,\cdot,\mu * \rho_n)* \rho_n(x) \quad \text{and} \quad \partial^2_x u^n(t,x,\mu) = \partial^2_x u(t,\cdot,\mu * \rho_n) * \rho_n(x).$$ These two functions are continuous on $[0,T]\times \R^d \times \PPP_2(\R^d)$ by the dominated convergence theorem and the fact that $u$ is jointly continuous. We define $\tilde{\rho_n}$ by $\tilde{\rho_n}(x,v) := \rho_n(x) \rho_n(v)$ for all $x,v \in \R^d.$ It is easy to see that $(\tilde{\rho_n})_n$  is a mollifying sequence on $\R^{2d}.$ Next, we claim that for all $(t,x) \in [0,T] \times \R^d$, $u^n(t,x,\cdot)$ has a linear derivative given by  \begin{equation}\label{eqderlin}
	\frac{\delta u^n}{\delta m} (t,x,\mu)(v) := \del (t,\cdot, \mu * \rho_n)(\cdot) *\tilde{\rho_n}(x,v).\end{equation}  This convolution is well-defined as $\del$ is jointly continuous. To prove \eqref{eqderlin}, note first that the bound of Assumption $(3)$ in Definition $\ref{sobolevW}$ implies that for all $(t,x) \in [0,T] \times \R^d$, $\frac{\delta u^n}{\delta m} (t,x,\mu)(\cdot)$ is at most of quadratic growth, uniformly in $\mu$ on each compact set.  Since for all $(t,x) \in [0,T] \times \R^d$, $\del(t,x,\cdot)(\cdot)$ is continuous on $\PPP_2(\R^d) \times \R^d$, the dominated convergence theorem proves that $\frac{\delta u^n}{\delta m} (t,x,\cdot)(\cdot)$ is  continuous. As explained in Remark \ref{Rqlinearderivative}, it is enough to compute, for $\mu,\nu \in \PPP_2(\R^d)$ and $\lambda \in [0,1]$, the derivative with respect to $\lambda$ of $u^n(t,x,m_{\lambda})$, where $m_{\lambda} = \lambda \mu + (1- \lambda) \nu.$ As recalled in the proof of Theorem \ref{ito_krylov}, when $(t,x)$ are fixed  $$\frac{d }{d\lambda} u(t,x,m_{\lambda} * \rho_n) = \int_{\R^d}\del(t,x,m_{\lambda}* \rho_n)*\rho_n(v) \, d(\mu - \nu)(v).$$ Thanks to the bound Assumption $(3)$ in Definition $\ref{sobolevW}$ for all compact $K \subset \R^d,$ one has  $$ \sup_{x \in K}\sup_{\lambda \in [0,1]} \left| \frac{d }{d\lambda} u(t,x,m_{\lambda} * \rho_n)\right| \leq C \left( 1 + \int_{\R^d} |v|^2 \, d(\mu + \nu)(v)\right).$$We can conclude with the help of Leibniz's rule and Fubini's theorem that   $$\frac{d }{d\lambda} u^n(t,x,m_{\lambda}) = \int_{\R^d}\del(t,\cdot,m_{\lambda}* \rho_n)*\tilde{\rho}_n(x,v) \, d(\mu - \nu)(v).$$ It follows from the joint continuity of $\del$ and Leibniz's rule that $\frac{\delta u^n}{\delta m}$ is $\CC^2$ with respect to $v$ and that  $$\left\{ \begin{aligned}
	\partial_v \frac{\delta u^n}{\delta m}(t,x,\mu)(v) &=& \frac{\delta u}{\delta m}(t,\cdot,\mu*\rho_n)(\cdot) * \partial_v\tilde{\rho_n}(x,v) = \partial_v\frac{\delta u}{\delta m}(t,\cdot,\mu*\rho_n)(\cdot) * \tilde{\rho_n}(x,v) \\\vspace{10pt} \partial_v^2 \frac{\delta u^n}{\delta m}(t,x,\mu)(v) &=& \frac{\delta u}{\delta m}(t,\cdot,\mu*\rho_n)(\cdot) * \partial^2_v\tilde{\rho_n}(x,v) = \partial_v^2\frac{\delta u}{\delta m}(t,\cdot,\mu*\rho_n)(\cdot) * \tilde{\rho_n}(x,v). \end{aligned}\right.$$ 
	Note that $\partial_v \frac{\delta u^n}{\delta m}$ and $\partial^2_v \frac{\delta u^n}{\delta m}$ are continuous on $[0,T] \times \R^d \times \PPP_2(\R^d) \times \R^d$ thanks to the dominated convergence theorem and the joint continuity of $\del.$ Moreover for all compact $ \KK \subset \PPP_2(\R^d)$ and for all $M>0$   \begin{equation}\label{eqbound}
		 \sup_{t \in[0,T]} \sup_{\mu \in \KK} \sup_{|x| \leq M}\sup_{v \in \R^d} \left|  \partial_v \frac{\delta u^n}{\delta m} (t,x,\mu)(v)\right| + \left|  \partial^2_v\frac{\delta u^n}{\delta m} (t,x,\mu)(v)\right| < + \infty. \end{equation} Indeed,  Hölder's inequality ensures that  \begin{align*}
		 &\sup_{|x| \leq M}\sup_{v \in \R^d} \left|  \partial_v \frac{\delta u^n}{\delta m} (t,x,\mu)(v)\right| + \left|  \partial^2_v\frac{\delta u^n}{\delta m} (t,x,\mu)(v)\right| \\&\leq \left[ \left\Vert \partial_v\frac{\delta u}{\delta m}(t,\cdot,\mu*\rho_n)(\cdot) \right\Vert_{L^{k_2}(B_{M+1}\times \R^d)} + \left\Vert \partial_v^2\frac{\delta u}{\delta m}(t,\cdot,\mu*\rho_n)(\cdot) \right\Vert_{L^{k_2}(B_{M+1}\times \R^d)} \right] \Vert \tilde{\rho_n}\Vert_{L^{k_2'}(\R^{2d})},\end{align*} the ball $B_{M+1}$ coming from the fact that the support of $\tilde{\rho_n}$ is included in $B_1.$ Since $\KK * \rho_n$ is compact in $\PPP_2(\R^d)$ and included in $\PP(\R^d)$, Assumption $(5)$ in Definition \ref{sobolevW} ensures that there exists $C>0$ such that for all $\mu \in \KK$   $$ \begin{aligned}
		\sup_{t \in [0,T]}\sup_{|x| \leq M}\sup_{v \in \R^d} \left|  \partial_v \frac{\delta u^n}{\delta m} (t,x,\mu)(v)\right| + \left|  \partial^2_v\frac{\delta u^n}{\delta m} (t,x,\mu)(v)\right| &\leq C \left( 1 + \left\Vert \frac{d \mu * \rho_n}{dx} \right\Vert_{L^{k_2'}(\R^d)}^{\alpha_2}\right)  \Vert \tilde{\rho_n}\Vert_{L^{k_2'}(\R^{2d})}.
		 \end{aligned}$$
		 But we know that  $ \frac{d \mu * \rho_n}{dx} (x) = \displaystyle\int_{\R^d} \rho_n (x-y) \, d\mu(y).$ We conclude with Jensen's inequality that $$ \left\Vert \frac{d \mu * \rho_n}{dx} \right\Vert_{L^{k_2'}(\R^d)}^{\alpha_2} \leq \Vert \rho_n \Vert_{L^{k_2'}(\R^d)}^{\alpha_2}.$$ This proves \eqref{eqbound}.\\

	\noindent\textbf{Step 4: Itô's formula \eqref{formulaito2} for $\boldsymbol{u^n}$ when the coefficients $\boldsymbol{b}$ and $\boldsymbol{\sigma}$ are continuous. } \\
	
	We claim that $(t,x) \mapsto U^n(t,x) := u^n(t,x,\mu_t) \in \CC^{1,2}([0,T] \times \R^d)$. The regularity with respect to $x$ is clear with the preceding properties on $u^n$. Let us thus focus on the regularity with respect to the time variable. For $(t,x) \in [0,T]\times \R^d$ fixed, the regularity assumption on $u$ with respect to $t$ and the standard Itô formula for a flow of measures applied to $u^n(t,x,\cdot)$ (see Theorem 5.99 in Chapter $5$ of \cite{CarmonaProbabilisticTheoryMean2018}) ensure that we have for $h \in \R$ satisfying $t+h \geq 0$   
	\begin{align}\label{formulaitou_n2} 
	\notag u^n(t+h,x,\mu_{t+h}) - u^n(t,x,\mu_t)&=  u^n(t+h,x,\mu_{t+h}) - u^n(t,x,\mu_{t+h}) + u^n(t,x,\mu_{t+h}) - u^n(t,x,\mu_t) \\ &=\int_t^{t+h} \partial_t u^n(s,x,\mu_{t+h}) \, ds + \int_t^{t+h} \E \left( \partial_v \frac{\delta u^n}{\delta m} (t,x,\mu_s)(X_s)\cdot b_s\right) \, ds\\ \notag&\quad+ \frac12 \int_t^{t+h} \E \left( \partial^2_v \frac{\delta u^n}{\delta m} (t,x,\mu_s)(X_s)\cdot a_s\right) \, ds.
	\end{align}
	The function $(s,x,\mu) \in [0,T] \times \R^d \times \PPP_2(\R^d) \mapsto \partial_t u^n(s,x,\mu)$ is continuous so  $$ \frac1h \int_t^{t+h} \partial_t u^n(s,x,\mu_{t+h}) \, ds \underset{h \rightarrow 0}{\longrightarrow}  \partial_t u^n(t,x,\mu_{t}).$$ The two other terms in \eqref{formulaitou_n2} can be dealt similarly. Indeed, the dominated convergence theorem justified by  \eqref{eqbound} ensures that the functions $(s,x) \in [0,T] \times \R^d \mapsto  \E \left( \partial_v \frac{\delta u^n}{\delta m} (s,x,\mu_s)(X_s)\cdot b_s\right) $ and  $(s,x) \in [0,T] \times \R^d \mapsto  \E \left( \partial^2_v \frac{\delta u^n}{\delta m} (s,x,\mu_s)(X_s)\cdot a_s\right) $ are continuous. Then, it follows that $U^n \in \CC^{1,2}([0,T] \times \R^d)$ and that for all $(t,x) \in  [0,T] \times \R^d $  $$ \partial_t U^n(t,x) = \partial_t u^n(t,x,\mu_t)+ \E \left( \partial_v \frac{\delta u^n}{\delta m} (t,x,\mu_t)(X_t)\cdot b_t\right) + \frac12  \E \left( \partial^2_v \frac{\delta u^n}{\delta m} (t,x,\mu_t)(X_t)\cdot a_t\right).$$ We can now apply the classical Itô formula for $U^n$ and $\xi$, up to the random time $\tau$ defined at the end of Step 2, to obtain that almost surely, for all $t \in [0,T]$  \begin{align}\label{eqstep4}
	\notag u^n(t\wedge \tau,\xi_{t\wedge \tau},\mu_{t\wedge \tau}) &= u^n(0,\xi_0, \mu_0) + \int_0^{t\wedge \tau} \partial_t u^n(s,\xi_s,\mu_s) + \partial_x u^n(s,\xi_s,\mu_s)\cdot\eta_s + \frac12 \partial^2_x u^n(s,\xi_s,\mu_s) \cdot \gamma_s\gamma_s^* \, ds \\ &\quad+ \int_0^{t\wedge \tau} \tilde{\E} \left(\partial_v \frac{\delta u^n}{\delta m} (s,\xi_s,\mu_s)(\tilde{X_s})\cdot \tilde{b_s}\right) \, ds + \frac12 \int_0^{t\wedge \tau} \tilde{\E} \left(\partial^2_v \frac{\delta u^n}{\delta m} (s,\xi_s,\mu_s)(\tilde{X_s})\cdot  \tilde{a_s}\right) \, ds \\ \notag&\quad+ \int_0^{t\wedge \tau} \partial_x u^n(s,\xi_s,\mu_s)\cdot(\gamma_s \, dB_s). 
	\end{align}
	Note that \eqref{eqstep4} does not require Assumptions \textbf{(A)} and \textbf{(B)} on the Itô process $X.$ These assumptions will only be used in Step 6.\\

	\noindent\textbf{Step 5: Removing the continuity hypothesis on the coefficients $\boldsymbol{b}$ and $\boldsymbol{\sigma}$.} \\
	
	We consider $(b^m)_m$ and $(\sigma^m)_m$ two sequences of continuous and progressively measurable processes such that $$ \E \int_0^T |b^n_s - b_s|^2 + |\sigma^n_s - \sigma_s|^4 \, ds \rightarrow 0.$$ We set, for $t \leq T, $  $X^m_t := X_0 + \int_0^t b^m_s \, ds + \int_0^t \sigma^m_s \, dB_s,$ and $\mu^m_t$ the law of $X^m_t$. Owing to Step $4$, Itô's formula \eqref{eqstep4} holds true for $X^m$ and $\xi$. Now, we aim at taking the limit $m \rightarrow +\infty$ in \eqref{eqstep4}. Note that the set $\KK:=\{\mu^m_s, \, s\leq T, \, m \geq 1 \} \cup \{ \mu_s, \, s \leq T\}$ is compact in $ \PPP_2(\R^d)$. Indeed, using Jensen's inequality and the Burkholder-Davis-Gundy (BDG) inequalities, it is clear that  $\E \, \underset{t \leq T}{\sup} |X^m_t - X_t|^2 \rightarrow 0,$ thus $\underset{t \leq T}{\sup} \,W_2(\mu^m_t,\mu_t) \rightarrow 0$. We deduce that almost surely, for all $t \in [0,T] $  $$u^n(t,\xi_t,\mu^m_t) \underset{m \rightarrow + \infty}{\longrightarrow} u^n(t,\xi_t,\mu_t).$$ Now, we take the limit $m \rightarrow + \infty$ in the integrals in Itô's formula \eqref{eqstep4}.\\
	
	\noindent\textbf{(i) Martingale term in \eqref{eqstep4}.} Using BDG's inequality, there exists $C>0$ such that  \begin{align*}
	&\E \sup_{t \leq T} \left| \int_0^{t\wedge \tau} (\partial_x u^n(s,\xi_s,\mu^m_s) - \partial_x u^n(s,\xi_s,\mu_s))\cdot(\gamma_s \, dB_s)\right|^2\\ &\leq C \E \int_0^{T\wedge \tau}|\partial_x u^n(s,\xi_s,\mu^m_s) - \partial_x u^n(s,\xi_s,\mu_s)|^2 |\gamma_s|^2 \, ds \\ &\leq C \E \int_0^{T}|\partial_x u^n(s,\xi_s,\mu^m_s) - \partial_x u^n(s,\xi_s,\mu_s)|^2\1_{B_M}(\xi_s) |\gamma_s|^2 \, ds. 
	\end{align*}
	The dominated convergence theorem can be applied since $\gamma$ is bounded and $\partial_x u^n$ is jointly continuous on $[0,T]\times \R^d \times \PPP_2(\R^d).$ It shows that, up to an extraction, almost surely  $$ \forall t \leq T, \, \int_0^{t\wedge \tau} \partial_x u^n(s,\xi_s,\mu^m_s)\cdot (\gamma_s \, dB_s) \underset{m \rightarrow + \infty}{\longrightarrow} \int_0^{t\wedge \tau} \partial_x u^n(s,\xi_s,\mu_s)\cdot (\gamma_s \, dB_s).$$

	\noindent \textbf{(ii) Terms involving the linear derivative in \eqref{eqstep4}.} Let us write  \begin{align*}
	&\left|\int_0^{t\wedge \tau} \tilde{\E} \left( \partial_v \frac{\delta u^n}{\delta m}(s,\xi_s,\mu^m_s)(\tilde{X^m_s})\cdot \tilde{b_s^m} \right) \, ds - \int_0^{t\wedge \tau} \tilde{\E} \left( \partial_v \frac{\delta u^n}{\delta m}(s,\xi_s,\mu_s)(\tilde{X_s})\cdot  \tilde{b_s} \right) \, ds \right| \\ &\leq \int_0^{T\wedge \tau} \tilde{\E} \left| \partial_v \frac{\delta u^n}{\delta m}(s,\xi_s,\mu^m_s)(\tilde{X_s^m})\right||\tilde{b_s^m} -\tilde{b_s}| \, ds \\ &\quad+  \int_0^{T\wedge \tau} \tilde{\E} \left|\partial_v \frac{\delta u^n}{\delta m}(s,\xi_s,\mu^m_s)(\tilde{X_s^m}) - \partial_v \frac{\delta u^n}{\delta m}(s,\xi_s,\mu_s)(\tilde{X_s}) \right||\tilde{b_s} | \, ds \\ &=: I_1+I_2
	\end{align*}
	Cauchy-Schwarz's inequality ensures that  \begin{align*}
	I_1 \leq \left( \int_0^{T\wedge \tau} \tilde{\E} \left|\partial_v \frac{\delta u^n}{\delta m}(s, \xi_s,\mu^m_s)(\tilde{X_s^m})\right|^2 \, ds\right)^{1/2}\left(\int_0^{T} \tilde{\E}|\tilde{b_s^m} - \tilde{b_s}|^2 \, ds\right)^{1/2}.
	\end{align*}
	We conclude that $I_1$ converges to $0$ thanks to the bound \eqref{eqbound} proved in Step $3$ and since $\xi$ is bounded by $M$ on the set $\{\tau >0 \}.$
	To show that $I_2 \rightarrow 0,$ we use the fact that $b$ is bounded by $K$ to get \begin{align*}
	I_2 \leq K\int_0^{T\wedge \tau} \tilde{\E} \left|  \partial_v \frac{\delta u^n}{\delta m}(s,\xi_s,\mu^m_s)(\tilde{X_s^m}) - \partial_v \frac{\delta u^n}{\delta m}(s,\xi_s,\mu_s)(\tilde{X_s})\right|\, ds.
	\end{align*}
	The continuity of $\partial_v \frac{\delta u^n}{\delta m}$ and the convergence in $L^2$ of $(\tilde{X^m_s})_m$ to $\tilde{X_s}$ ensure that for all $\omega \in \Omega,$ \newline $  \left|  \partial_v \frac{\delta u^n}{\delta m}(s,\xi_s(\omega),\mu^m_s)(\tilde{X_s^m}) - \partial_v \frac{\delta u^n}{\delta m}(s,\xi_s(\omega),\mu_s)(\tilde{X_s})\right|$ converges in probability on $\tilde{\Omega}$ to $0$ as $m$ goes to infinity. Using a uniform integrability argument coming from \eqref{eqbound}, we deduce that $I_2$ converges to $0.$  Following the same strategy, one has for all $t \in [0,T]$  
	$$ \int_0^{t\wedge \tau} \tilde{\E} \left( \partial^2_v \frac{\delta u^n}{\delta m}(s,\xi_s,\mu^m_s)(\tilde{X^m_s})\cdot \tilde{a_s^m} \right) \, ds \underset{m \rightarrow + \infty}{\longrightarrow} \int_0^{t\wedge \tau} \tilde{\E} \left( \partial^2_v \frac{\delta u^n}{\delta m}(s,\xi_s,\mu_s)(\tilde{X_s})\cdot\tilde{a_s} \right) \, ds.$$
	
		\noindent \textbf{(iii) Terms involving standard derivatives in \eqref{eqstep4}.} It follows from the dominated convergence theorem that almost surely, for all $t \leq T$  \begin{align*}&\int_0^{t\wedge \tau} ( \partial_t u^n(s,\xi_s,\mu^m_s) + \partial_x u^n(s,\xi_s,\mu^m_s)\cdot\eta_s) \, ds + \frac12 \int_0^{t\wedge \tau} \partial^2_x u^n(s,\xi_s,\mu^m_s)\cdot \gamma_s\gamma_s^* \, ds\\   &\underset{m \rightarrow + \infty}{\longrightarrow} \int_0^{t\wedge \tau} ( \partial_t u^n(s,\xi_s,\mu_s) + \partial_x u^n(s,\xi_s,\mu_s)\cdot \eta_s) \, ds + \frac12 \int_0^{t\wedge \tau} \partial^2_x u^n(s,\xi_s,\mu_s)\cdot \gamma_s\gamma_s^* \, ds. \end{align*} 
	Indeed the functions $ \partial_t u^n$, $\partial_x u^n$ and $\partial_x^2 u^n$ are jointly continuous on $[0,T]\times \R^d \times \PPP_2(\R^d)$ and thus uniformly bounded on $[0,T]\times B_M \times \{\mu_s^m, \, s \in [0,T],\, m \geq 1 \}.$ Moreover, $\eta$ and $\gamma$ are also uniformly bounded.\\
	
	This concludes Step $5$.\\

	\noindent\textbf{Step 6: Letting $\boldsymbol{n \rightarrow + \infty.}$}\\
	
	From Step $5$, we deduce that Itô's formula \eqref{eqstep4} in Step 4 holds for $u^n$ up to time $\tau.$ To conclude the proof, we need to take the limit $n \rightarrow + \infty$ in each term of \eqref{eqstep4}. Then it remains to remove the stopping time $\tau$ as explained at the end of Step $2$ (i.e. letting $\tau \rightarrow T$). The continuity of $u$ ensures that almost surely, for all $t \leq T$, $u^n(t,\xi_t,\mu_t) \rightarrow u(t,\xi_t,\mu_t).$ We now focus on the integrals in Itô's formula \eqref{eqstep4}.\\
	
	\noindent \textbf{(i) Martingale term in \eqref{eqstep4}.} Thanks to BDG's inequality, Hölder's inequality, and the boundedness of $\gamma,$ we have  \begin{align*}
	&\E \sup_{t \leq T} \left| \int_0^{t \wedge \tau} (\partial_x u^n(s,\xi_s,\mu_s) - \partial_x u(s,\xi_s,\mu_s))\cdot(\gamma_s \, dB_s)\right|^2  \\&\leq C  \E\int_0^{T } |\partial_x u(s,\cdot,\mu_s* \rho_n)*\rho_n(\xi_s) - \partial_x u(s,\xi_s,\mu_s)|^2\1_{B_M}(\xi_s) \, ds\\ &= C \int_0^T\int_{B_M} |\partial_{x} u(s,\cdot,\mu_s* \rho_n)*\rho_n(x) -\partial_{x} u(s,x,\mu_s)|^2 q(s,x) \, dx \,ds \\ &\leq C\int_0^T \Vert \partial_{x} u(s,\cdot,\mu_s* \rho_n)*\rho_n -\partial_{x} u(s,\cdot,\mu_s) \Vert_{L^{2k_1}(B_M)}^2 \Vert q(s,\cdot) \Vert_{L^{k_1'}(B_M)} \, ds\\ &\leq C\int_0^T \Vert \partial_{x} u(s,\cdot,\mu_s* \rho_n)*\rho_n -\partial_{x} u(s,\cdot,\mu_s)* \rho_n \Vert_{L^{2k_1}(B_M)}^2 \Vert q(s,\cdot) \Vert_{L^{k_1'}(B_M)}  \, ds \\   &\quad+ C\int_0^T \Vert \partial_{x} u(s,\cdot,\mu_s)*\rho_n -\partial_{x} u(s,\cdot,\mu_s) \Vert_{L^{2k_1}(B_M)}^2 \Vert q(s,\cdot) \Vert_{L^{k_1'}(B_M)} \, ds \\ &=: I_1 + I_2.
	\end{align*}
	We prove that $I_1$ and $I_2$ converge to $0.$ First note that, due to the convolution inequality $L^r*L^1$, we have for $f \in L^r_{\text{loc}}(\R^d)$ and for all $R>0,$ $\|f * \rho_n \|_{L^r(B_R)} \leq \|f\|_{L^r(B_{R+1})}.$ The control on $B_{R+1}$ follows from the fact that the support of each $\rho_n$ is included in $B_1$. Hence
	
	\begin{align*}
I_1 & \leq C \int_0^T  \Vert \partial_{x} u(s,\cdot,\mu_s* \rho_n) -\partial_{x} u(s,\cdot,\mu_s) \Vert_{L^{2k_1}(B_{M+1})}^2 \Vert q(s,\cdot) \Vert_{L^{k_1'}(B_{M+1})}  \, ds \, =: \tilde{I_1}.
	\end{align*} As a consequence of Sobolev embedding theorem, for all $t$, the function  $$ \mu \in  (\PP(\R^d),d_{\PP}) \mapsto \partial_x u(t,\cdot,\mu) \in L^{\infty}(B_{M+1})$$ is continuous. Since $\mu_s \in \PP(\R^d)$ for almost all $s$ and thanks to Assumption $(2)$ in Definition \ref{defespaceP}, we deduce that the integrand in $\tilde{I_1}$ converges to $0$ for almost all $s$. It follows from the dominated convergence theorem (see \eqref{bound1} in Step 1) that $\tilde{I_1}$ converges to $0$, as well as $I_1.$  We now focus on $I_2.$ The integrand in $I_2$ converges to $0$ for almost all $s$ because $\partial_x u (s,\cdot,\mu_s) \in L^{2k_1}(B_{M}).$ We conclude with the dominated convergence theorem as previously. This shows that, up to an extraction, almost surely  $$ \sup_{t \leq T} \left| \int_0^{t \wedge \tau} \partial_x u^n(s,\xi_s,\mu_s) \cdot(\gamma_s \, dB_s) - \int_0^{t \wedge \tau} \partial_x u(s,\xi_s,\mu_s) \cdot(\gamma_s \, dB_s) \right| \rightarrow 0. $$ 
	
	\noindent \textbf{(ii) Terms involving the linear derivative in \eqref{eqstep4}.} Following the same strategy, we obtain using Hölder's inequality \begin{align*}
	&\E \, \sup_{t \leq T} \left|\int_0^{t \wedge \tau}\tilde{\E} \left(\partial_v \frac{\delta u^n}{\delta m} (s,\xi_s,\mu_s)(\tilde{X_s})\cdot\tilde{b_s}\right)  \, ds - \int_0^{t \wedge \tau}\tilde{\E} \left(\partial_v \frac{\delta u}{\delta m} (s,\xi_s,\mu_s)(\tilde{X_s})\cdot\tilde{b_s}\right)  \, ds \right|\\ &\leq \E\tilde{\E}\int_0^{T \wedge \tau}\left|\partial_v \frac{\delta u^n}{\delta m} (s,\xi_s,\mu_s)(\tilde{X_s}).\tilde{b_s} -\partial_v \frac{\delta u}{\delta m} (s,\xi_s,\mu_s)(\tilde{X_s}).\tilde{b_s}\right| \, ds \\ &\leq \E\tilde{\E}\int_0^{T }\left|\partial_v \frac{\delta u^n}{\delta m} (s,\xi_s,\mu_s)(\tilde{X_s})\cdot\tilde{b_s} -\partial_v \frac{\delta u}{\delta m} (s,\xi_s,\mu_s)(\tilde{X_s})\cdot\tilde{b_s}\right|\1_{B_M}(\xi_s) \, ds\\ &\leq C \int_0^T \int_{B_M\times \R^d}  \left|\partial_v \frac{\delta u^n}{\delta m} (s,x,\mu_s)(v) -\partial_v \frac{\delta u}{\delta m} (s,x,\mu_s)(v)\right| q(s,x)p(s,v) \, dx \, dv \, ds \\ &\leq C \int_0^T \left\Vert  \partial_v \frac{\delta u}{\delta m} (s,\cdot,\mu_s* \rho_n)(\cdot)* \tilde{\rho_n} -\partial_v \frac{\delta u}{\delta m} (s,\cdot,\mu_s)(\cdot) \right\Vert_{L^{k_2}(B_M\times\R^d)} \Vert q(s,\cdot) \Vert_{L^{k_2'}(B_M)}  \Vert p(s,\cdot) \Vert_{L^{k_2'}(\R^d)}  \, ds.
	\end{align*}
	 The dominated convergence theorem justified by Assumption $(4)$ in Definition \ref{sobolevW} and \eqref{bound1} in Step $1$ ensures that this term converges to $0.$ The same argument holds true for the term involving $\partial^2_v \del$. \\
	
	\noindent \textbf{(iii) Terms involving standard derivatives in \eqref{eqstep4}.} The convergence of the term involving $\partial_tu^n$ in \eqref{eqstep4} follows from the continuity of $\partial_t u$ on $[0,T] \times \R^d \times \PPP_2(\R^d)$ and the dominated convergence theorem since almost surely on the set $\{ \tau >0 \}$  $$ \sup_{s \in [0,T]}\sup_{n\geq1} |\partial_t u^n(s,\xi_s,\mu_s) | \leq  \sup_{s \in [0,T]}\sup_{n\geq 1} \sup_{|x| \leq M+1} |\partial_t u(s,x,\mu_s* \rho_n) | <+ \infty.$$ For the spatial derivatives, Hölder's inequality ensures that \begin{align*}
	&\E \, \sup_{t \leq T} \left| \int_0^{t \wedge \tau} \partial_x u^n(s,\xi_s,\mu_s)\cdot\eta_s \, ds - \int_0^{t \wedge \tau} \partial_x u(s,\xi_s,\mu_s)\cdot\eta_s \, ds \right|  \\ &\leq C\, \int_0^{T} \Vert  \partial_x u(s,\cdot,\mu_s* \rho_n)*\rho_n-\partial_x u(s,\cdot,\mu_s)\Vert_{L^{k_1}(B_M)} \Vert q(s,\cdot)\Vert_{L^{k_1'}(B_M)} \, ds.
	\end{align*}
	The right-hand side term converges to $0$ with same reasoning as before. This shows that, up to an extraction, one has almost surely  $$ \sup_{t \leq T} \left| \int_0^{t \wedge \tau} \partial_x u^n(s,\xi_s,\mu_s)\cdot\eta_s \, ds - \int_0^{t \wedge \tau}\partial_x u(s,\xi_s,\mu_s)\cdot\eta_s \, ds \right|  \underset{n \rightarrow + \infty}{\longrightarrow} 0.$$ The term involving $\partial^2_x u$ in \eqref{eqstep4} is dealt similarly.\\
	
	 Taking the limit $ n \rightarrow + \infty$ in  \eqref{eqstep4}, up to an extraction,  we conclude that almost surely, for all $t \in [0,T]$  \begin{align*}
	u(t\wedge \tau,\xi_{t\wedge \tau},\mu_{t\wedge \tau}) &= u(0,\xi_0, \mu_0) \\ &\quad+ \int_0^{t\wedge \tau} ( \partial_t u(s,\xi_s,\mu_s) + \partial_x u(s,\xi_s,\mu_s)\cdot\eta_s) \, ds + \frac12 \int_0^{t\wedge \tau} \partial^2_x u(s,\xi_s,\mu_s) \cdot \gamma_s\gamma_s^* \, ds \\ &\quad+ \int_0^{t\wedge \tau} \tilde{\E} \left(\partial_v \frac{\delta u}{\delta m} (s,\xi_s,\mu_s)(\tilde{X_s})\cdot\tilde{b_s}\right) \, ds + \frac12 \int_0^{t\wedge \tau} \tilde{\E} \left(\partial^2_v \frac{\delta u}{\delta m} (s,\xi_s,\mu_s)(\tilde{X_s})\cdot \tilde{a_s}\right) \, ds \\ &\quad+ \int_0^{t\wedge \tau} \partial_x u(s,\xi_s,\mu_s)\cdot(\gamma_s \, dB_s).
	\end{align*} This ends the proof as explained in Step 2. \hfill$\square$

	\appendix
	\section{}\label{sectionappendix}
	
\subsection{Proof of Example \ref{choicedistance}}\label{proofchoicedistance}

\textbf{(1)} It follows from the contraction inequality in Lemma \ref{convolutionmesure} and Corollary \ref{regularisationconvolutionmeasure}.
	
	\textbf{(2)} To prove (\textbf{H1}), we fix $n \geq 1$ and $\mu_j \overset{W_2}{\longrightarrow} \mu \in \PPP_2(\R^d).$ For $\nu \in \PPP_2(\R^d),$ the density of $\nu* \rho_n$ is given by  $$x \in \R^d \mapsto \rho_n * \nu(x) = \int_{\R^d}\rho_n(x-y) \, d\nu(y).$$ Hence,  $$d_k(\mu_j* \rho_n,\mu* \rho_n) = \left\Vert  \int_{\R^d}\rho_n(\cdot-y) \, d\mu_j(y) -  \int_{\R^d}\rho_n(\cdot-y) \, d\mu(y) \right\Vert_{L^{k'}(\R^d)}.$$ Using Lemma \ref{continuiteconvolution}, we conclude that  $d_k(\mu_j* \rho_n,\mu* \rho_n) \underset{j \rightarrow + \infty}{\longrightarrow}0.$
	 For \textbf{(H2)}, let $\mu \in \PP(\R^d)$ and denote by $f\in L^{k'}(\R^d)$ the density of $\mu.$ For $n\geq 1$, we have  $$ \frac{d\mu*\rho_n}{dx} = f* \rho_n \overset{L^{k'}}{\longrightarrow} f,$$ owing to Lemma \ref{regularisation}.
	 
	 \hfill$\square$	
	 
	 \subsection{Proof of Example \ref{exquadratic}}\label{proofexquadratic}
	
	 	Let us give the detailed proof in the bilinear case $N=2$. It is standard (see Example 4 page 389 in Chapter $5$ of \cite{CarmonaProbabilisticTheoryMean2018}) that $u$ has a linear derivative given by  $$ \del (\mu)(v) = \int_{\R^d} g(v,y) \, d\mu(y) + \int_{\R^d} g(y,v) \, d\mu(y).$$ We will only treat the first term since the other one can be dealt similarly.\\
	 	
	 	\noindent\textbf{Computation of the distributional derivatives and continuity:} Let $\mu \in \PP(\R^d)$ and $f \in L^{(d+1)'}(\R^d)$ be its density. By interpolation, we know that $f \in L^{r'}(\R^d)$ for all $r \geq d+1.$ Let $\varphi \in \CC^{\infty}_c(\R^d)$ and $i \in \{1, \dots, d\}.$ Using Fubini's theorem, justified  by the quadratic growth of $g$ and the fact that $f\,dx \in \PPP_2(\R^d),$ we have  $$ \int_{\R^d} \left( \int_{\R^d} g(x,y) f(y) \, dy\right) \partial_{v_i} \varphi(v) \, dv =  \int_{\R^d \times \R^d}g(v,y)f(y) \partial_{v_i} \varphi(v) \, dy \, dv.$$ Let us define $f_n(x) = \frac{1}{\mu(B_n)} (f\1_{B_n})*\rho_n(x),$ for $n$ large enough to have $\mu(B_n) >0.$ The function $f_n$ is a probability density which is in $\CC^{\infty}_c(\R^d).$ It easily follows from  Lemma \ref{convolution}, Lemma \ref{regularisation} and the dominated convergence theorem that \begin{equation}\label{eqapprox}
	 		 f_n \overset{L^{k'}}{\longrightarrow} f \quad \text{and} \quad f_n \overset{W_2}{\longrightarrow} f. \end{equation} For a fixed $n \geq 1$, we have by definition of the distributional derivative \begin{equation}\label{eqcomputationderivative}
	 		 	 \int_{\R^d \times \R^d}g(v,y)f_n(y) \partial_{v_i} \varphi(v) \, dy \, dv = - \int_{\R^d \times \R^d}\partial_{v_i}g(v,y)f_n(y)  \varphi(v) \, dy \, dv. \end{equation}  Our aim is to take the limit $ n \rightarrow + \infty$ in both side of the previous equality. Using Fubini's theorem, the left-hand side term is equal to $$ \int_{\R^d} \left( \int_{\R^d} g(v,y) \partial_{v_i} \varphi(v) \, dv \right) f_n(y)\, dy.$$  Moreover, it converges to  $$ \int_{\R^d\times \R^d}  g(v,y) \partial_{v_i} \varphi(v)  f(y) \, dy \, dv.$$  Indeed, $f_n \overset{W_2}{\longrightarrow} f$ and the function $ y \mapsto \int_{\R^d} g(v,y) \partial_{v_i} \varphi(v) \, dv$ is continuous and at most of quadratic growth. For the right-hand side term, we prove that  $$ \int_{\R^d \times \R^d}\partial_{v_i}g(v,y)f_n(y)  \varphi(v) \, dy \, dv \rightarrow \int_{\R^d \times \R^d}\partial_{v_i}g(v,y)f(y)  \varphi(v) \, dy \, dv.$$ Note that the limit is well-defined using Hölder's inequality  \begin{align*}
	 	\int_{\R^d \times \R^d}|\partial_{v_i}g(v,y)f(y)  \varphi(v) |\, dy \, dv &\leq \Vert f \Vert_{L^{k'}(\R^d)} \int_{\R^d} \varphi(v) \Vert \partial_{v_i} g(v,\cdot)\Vert_{L^k(\R^d)} \, dv.
	 	\end{align*}
	 The right-hand side term is finite because $v \mapsto  \Vert \partial_{v_i} g(v,\cdot)\Vert_{L^k(\R^d)} \in L^k(\R^d).$ The same inequality shows that  $$ \left|	\int_{\R^d \times \R^d}\partial_{v_i}g(v,y)(f_n(y)-f(y) ) \varphi(v) \, dy \, dv \right| \leq \Vert f_n-f \Vert_{L^{k'}(\R^d)} \int_{\R^d} \varphi(v) \Vert \partial_{v_i} g(v,\cdot)\Vert_{L^k(\R^d)} \, dv \underset{n \rightarrow + \infty}{\longrightarrow} 0,$$ thanks to \eqref{eqapprox}. Taking the limit $n \rightarrow + \infty$ in \eqref{eqcomputationderivative}, we deduce that : $$ \int_{\R^d \times \R^d}g(v,y)f(y) \partial_{v_i} \varphi(v) \, dy \, dv = - \int_{\R^d \times \R^d}\partial_{v_i}g(v,y)f(y)  \varphi(v) \, dy \, dv. $$ Hence, the distributional derivative of  $v \mapsto \int_{\R^d} g(v,y) f(y) \, dy $ is given by the function $$ v \mapsto \int_{\R^d} \partial_v g(v,y) f(y) \, dy.$$ Moreover, it belongs to $L^k(\R^d)$ because applying Hölder's inequality, one has \begin{align*}
	 	\int_{\R^d} \left|\int_{\R^d} \partial_v g(v,y) f(y) \, dy\right|^k \, dv &\leq \int_{\R^d} \Vert \partial_v g(v,\cdot) \Vert_{L^{k}(\R^d)}^k \Vert f \Vert_{L^{k'}(\R^d)}^k \, dv \\ &= \Vert \partial_v g\Vert_{L^k(\R^d \times \R^d)}^k \Vert f \Vert_{L^{k'}(\R^d)}^{k}.
	 	\end{align*}
	 	Note that this inequality and the linearity in $f$ justify that  $ \mu \in (\PP(\R^d),d_k) \mapsto \int_{\R^d} \partial_v g(\cdot,y) \, d\mu(y) \in L^{k}(\R^d)$ is continuous with  \begin{equation}\label{eqcontinuity1}
	 		\left\Vert \partial_v \del(\mu) \right\Vert_{L^k(\R^d)} \leq  \left\Vert\frac{d\mu}{dx}\right\Vert_{L^{k'}(\R^d)} \Vert \nabla g \Vert_{L^k(\R^d\times \R^d)}.
	 	\end{equation} Following the same lines, we show that the distributional derivative of order $2$ of $\del (\mu)$, for $\mu \in \PP(\R^d)$, is given by the $\R^{d\times d}$-valued function  $$ v \mapsto \int_{\R^d} \partial^2_{v} g(v,y) \, d\mu(y) + \int_{\R^d} \partial^2_{y} g(y,v) \, d\mu(y).$$ It is also a continuous function from $(\PP(\R^d),d_k)$ into $ L^{k}(\R^d).$ Indeed, as previously, we obtain : \begin{equation}\label{eqcontinuity2}
	 		 \left\Vert \partial^2_v \del(\mu) \right\Vert_{L^k(\R^d)} \leq  \left\Vert\frac{d\mu}{dx}\right\Vert_{L^{k'}(\R^d)} \Vert \nabla^2 g \Vert_{L^k(\R^d\times \R^d)}.
	 	\end{equation}  
	 	
	 	\noindent\textbf{Growth property: } Using the inequalities \eqref{eqcontinuity1} and \eqref{eqcontinuity2} of the previous step, one has for all $\mu \in \PP(\R^d)$  \begin{align*} \left\Vert \partial_v \del(\mu) \right\Vert_{L^k(\R^d)} + \left\Vert \partial^2_v \del(\mu) \right\Vert_{L^k(\R^d)} &\leq  \left\Vert\frac{d\mu}{dx}\right\Vert_{L^{k'}(\R^d)} \left[\Vert \nabla g \Vert_{L^k(\R^d\times \R^d)} + \Vert \nabla^2 g \Vert_{L^k(\R^d\times \R^d)} \right].
	 	\end{align*}
	 	The second point in Definition \ref{sobolevW1} is thus satisfied with $\alpha =1$ because we have supposed that $k \geq 2d.$\\
	 	
	 	In the general case $N \geq 2$, one can show following the same lines that $u$ admits a linear derivative and that for all $\mu \in \PP(\R^d)$, its distributional derivative is given for all $v \in \R^d$ by $$ \partial_v \del (\mu)(v) = \sum_{j=1}^N \int_{(\R^d)^{N-1}} \partial_{x_j} g(x_1,\dots,x_{j-1},v,x_{j+1}, \dots,x_N) \, d\mu(x_1) \dots \, d\mu(x_{j-1}) \, d\mu(x_{j+1}) \dots \, d\mu(x_N).$$ 
	 	Denoting by $f$ the density of $\mu$ and using Hölder's inequality, we obtain as previously that for all $j \in \{1, \dots,N\}$
	 	\begin{align*}
	 	&\int_{\R^d} \left| \int_{(\R^d)^{N-1}} \partial_{x_j} g(x_1,\dots,x_{j-1},v,x_{j+1}, \dots,x_N) \, d\mu(x_1) \dots \, d\mu(x_{j-1}) \, d\mu(x_{j+1}) \dots \, d\mu(x_N)\right|^k \, dv\\  &= \Vert \partial_{x_j}g\Vert_{L^k((\R^d)^N)}^k \Vert f \Vert_{L^{k'}(\R^d)}^{(N-1)k}.
	 	\end{align*}
	 	
	 	We easily show that $\mu \in (\PP(\R^d),d_k) \mapsto \partial_v \del(\mu) \in L^k(\R^d)$ is continuous and the same properties hold for the distributional derivative of order two. We deduce that $\mu \in \PP(\R^d)$  \begin{align*} \left\Vert \partial_v \del(\mu) \right\Vert_{L^k(\R^d)} + \left\Vert \partial^2_v \del(\mu) \right\Vert_{L^k(\R^d)} &\leq  \left\Vert\frac{d\mu}{dx}\right\Vert_{L^{k'}(\R^d)}^{N-1} \left[\Vert \nabla g \Vert_{L^k((\R^d)^N)} + \Vert \nabla^2 g \Vert_{L^k((\R^d)^N)} \right].
	 	\end{align*}
	 	The second point in Definition \ref{sobolevW1} is thus satisfied with $\alpha =N-1$ because we have supposed that $k \geq Nd.$
	 	
	 	\hfill$\square$

	 	\subsection{Proof of Example \ref{exconvol}}\label{proofexconvol}
	
	Note that $f * \mu$ and $u(\mu)$ are well-defined for $\mu \in \PPP_2(\R^d).$ Indeed, it follows from Sobolev embedding theorem (see Corollary 9.14 in \cite{BrezisFunctionalAnalysisSobolev2010}) that $f \in \CC^1(\R^d,\R)$ and $\partial_x f \in (L^{\infty}(\R^d))^d$. Thus $f$ is at most of linear growth. Since $f$ is continuous and at most of linear growth, it is easy to see that $u$ has a linear derivative given by  $$ \forall \mu \in \PPP_2(\R^d), \, \forall v \in \R^d, \, \del (\mu)(v) = f * \mu(v) + \tilde{f} * \mu(v),$$ where $\tilde{f}(x) = f(-x)$ (see Example 2 page 386 in Chapter 5 of \cite{CarmonaProbabilisticTheoryMean2018}). An easy computation based on Fubini's theorem shows that the distributional derivatives of order $1$ and $2$ of $\del(\mu)$ are given by  $$\forall i,j \in \{1,\dots,d\},\, \left\{ \begin{array}{cll}\partial_{v_i} \del (\mu) &= \partial_{v_i} f * \mu+ \partial_{x_i} \tilde{f} * \mu \\ \partial_{v_i \, v_j} \del (\mu)&= \partial_{v_i \, v_j} f * \mu + \partial_{v_i \, v_j} \tilde{f} * \mu, \end{array}\right.$$ as elements of $L^{k+1}(\R^d).$ These functions are continuous with respect to $\mu \in \PPP_2(\R^d)$ owing to  Lemma \ref{continuiteconvolution}. It remains to apply the first point in Remark \ref{rqW1} to conclude.
	
	\hfill$\square$
	
	\subsection{Proof of Example \ref{ex-non-lin}} \label{proofex-non-lin}
	
	The function $u$ is well-defined and continuous because $\nabla g \in L^{\infty}(\R^d)$ and is continuous thanks to Sobolev embedding theorem. Thus $g$ is at most of linear growth. It follows from the continuity of $\mu \in \PPP_2(\R^d) \mapsto \int_{\R^d} g \, d\mu $ that the function $u$ admits a linear derivative given by  $$ \forall (\mu,v) \in \PPP_2(\R^d)\times \R^d, \, \del (\mu)(v) = g(v) F'\left(\int_{\R^d}g \, d\mu \right).$$
	 We thus have in the sense of distributions  $$ \forall \mu \in \PPP_2(\R^d), \, \forall v \in \R^d, \, \partial_v \del(\mu)(v) = \nabla g (v)  F'\left(\int_{\R^d}g \, d\mu \right).$$  Moreover, the function  $$ \mu \in \PPP_2(\R^d) \mapsto \partial_v \del(\mu)(\cdot) \in L^{k}( \R^d)$$ is continuous because $F \in \CC^1(\R;\R)$ and $\nabla g \in L^k(\R^d)$. The same reasoning proves that  $$ \forall \mu \in \PPP_2(\R^d), \, \forall v \in \R^d, \,  \partial_v^2 \del(\mu)(v) = \nabla^2 g (v)  F' \left(\int_{\R^d}g \, d\mu \right),$$ and that the function   $$ \mu \in \PPP_2(\R^d) \mapsto \partial_v^2 \del(\mu)(\cdot) \in L^{k}(\R^d)$$ is continuous. We conclude that $u\in \mathcal{W}_1(\R^d)$ with Remark \ref{rqW1}.
	
	\hfill$\square$
	
	\subsection{Proof of Example \ref{exlinear2}}\label{proofexlinear2}
		The function $u$ is well-defined and continuous. Indeed, Sobolev embedding theorem implies that $\nabla g \in L^{\infty}(\R^d)$ and is continuous. Hence $g$ is at most of linear growth. Following the same method as in the proof of Example \ref{exquadratic}, we obtain that  $$ \forall \mu \in \PP(\R^d),\, \partial_x u(\cdot,\mu) = \int_{\R^d} \partial_x g(\cdot,y)\, d\mu(y).$$ Moreover $$ \forall \mu \in \PP(\R^d), \, \Vert \partial_x u ( \cdot, \mu) \Vert_{L^k(\R^d)} \leq \Vert \nabla g \Vert_{L^{k}(\R^{2d})} \left\Vert \frac{d\mu}{dx} \right\Vert_{L^{k'}(\R^d)}.$$ This yields the continuity of the function $$ \mu \in (\PP(\R^d),d_k) \mapsto \partial_x u(\cdot,\mu) \in L^{k}(\R^d).$$ Moreover, keeping the notations of Definition \ref{sobolevW}, Assumption $(2)$ is satisfied and setting $ \alpha_1 = 1,$ Assumption $(5)$ is satisfied because we have supposed $k \geq 5d.$ The same holds true for $\partial^2_x u.$ Since $g$ is continuous and at most of linear growth, the linear derivative of $u$ satisfies Assumption $(3)$ in Definition \ref{sobolevW} and is given, for all $x,v \in \R^d$ and for all $\mu \in \PP_2(\R^d),$ by  $$ \del (x,\mu)(v) = g(x,v).$$ As $\nabla g \in (W^{1,k}(\R^{2d}))^d,$ Assumption $(4)$ in Definition \ref{sobolevW} is satisfied, as well as the growth property in Assumption $(5)$ with $\alpha_2=0.$
		 
		\hfill$\square$

	\subsection{Proof of Example \ref{ex1}}\label{proofex1}

	As in \ref{proofexlinear2}, the function $u$ is well-defined and continuous because $\nabla g \in L^{\infty}(\R^d)$ and is continuous. Thus $g$ is at most of linear growth. It is clear with the assumption on $\nabla F$ that for all $\mu \in \PPP_2(\R^d)$, $u(\cdot,\mu) \in W^{2,k_1}_{\text{loc}}(\R^d).$ It follows from the continuity of $\mu \in \PPP_2(\R^d) \mapsto \int_{\R^d} g \, d\mu $ that the function $$ \mu \in \PPP_2(\R^d)\mapsto \partial_x u(\cdot,\mu) = \partial_x F \left(\cdot,\int_{\R^d} g \, d\mu\right) \in (W^{1,k_1}(B_R))^d ,$$ is also continuous for all $R>0.$ Moreover, it is easy to show with Remark \ref{Rqlinearderivative} that for all $x \in \R^d,$ $u(x,\cdot)$ admits a linear derivative given by  $$ \forall (\mu,v) \in \PPP_2(\R^d)\times \R^d, \, \del (x,\mu)(v) = g(v) \partial_y F\left(x,\int_{\R^d}g \, d\mu \right).$$ Assumption $(3)$ in Definition \ref{sobolevW} is clearly satisfied because $\partial_y F$ is continuous. Next, we compute the derivatives of $\del(\cdot,\mu)(\cdot)$ with respect to $v$ in the sense of distributions. For $\phi \in \CC^{\infty}_c(\R^{2d})$ and $\mu \in \PPP_2(\R^d)$, Fubini's theorem ensures that \begin{align*}
		\int_{\R^{2d}}  g(v) \partial_y F\left(x,\int_{\R^d}g \, d\mu \right) \partial_v \phi(x,v) \, dx \, dv &= \int_{\R^d} \left(\int_{\R^d}g(v)\partial_v \phi(x,v) \, dv\right) \partial_y F\left(x,\int_{\R^d}g \, d\mu \right)  \, dx \\ &= -  \int_{\R^d} \left(\int_{\R^d}\nabla g(v) \phi(x,v) \, dv\right) \partial_y F\left(x,\int_{\R^d}g \, d\mu \right)  \, dx \\ &= - \int_{\R^{2d}} \left(\nabla g (v)  \partial_y F\left(x,\int_{\R^d}g \, d\mu \right) \right) \, \phi(x,v) \, dx\, dv.
	\end{align*}
	This proves exactly that  $$ \forall \mu \in \PPP_2(\R^d), \, \forall x,v \in \R^d, \, \partial_v \del(x,\mu)(v) = \nabla g (v)  \partial_y F\left(x,\int_{\R^d}g \, d\mu \right).$$ Since $ \nabla g \in L^{k_2}(\R^d)$ and $\partial_y F\left(\cdot,\int_{\R^d}g \, d\mu \right) \in L^{\infty}(B_R),$ for all $R>0$ and $\mu \in \PPP_2(\R^d),$ the function  $$(x,v) \in B_R \times \R^d \mapsto \partial_v \del(x,\mu)(v) $$ belongs to $L^{k_2}(B_R\times \R^d).$  Moreover, the function  $$ \mu \in \PPP_2(\R^d) \mapsto \partial_v \del(\cdot,\mu)(\cdot) \in L^{k_2}(B_R\times \R^d)$$ is continuous because $F \in \CC^1(\R^d \times \R;\R)$ and thus $ y \mapsto \partial_y  F(\cdot, y) \in L^{\infty}(B_R)$ is continuous. The same reasoning proves that  $$ \forall \mu \in \PPP_2(\R^d), \, \forall x,v \in \R^d, \,  \partial_v^2 \del(x,\mu)(v) = \nabla^2 g (v)  \partial_y F\left(x,\int_{\R^d}g \, d\mu \right),$$ and that the function   $$ \mu \in \PPP_2(\R^d) \mapsto \partial_v^2 \del(\cdot,\mu)(\cdot) \in L^{k_2}(B_R\times \R^d)$$ is continuous for all $R>0.$ We conclude that $u\in \mathcal{W}_2(\R^d)$ with Remark \ref{CorW_2}.
	
	\hfill$\square$
	
	\section*{Acknowledgements}
	
	I would like to thank Paul-Eric Chaudru de Raynal and Mihai Gradinaru for their supervision, advices and also for their careful reading of the paper. I would also like to thank Emeline Luirard for her careful reading of the paper.

	\bibliographystyle{plain}
	\bibliography{bibli}

\end{document}